\documentclass[review]{elsarticle}
\usepackage[letterpaper,margin=1in]{geometry}

\usepackage[dvipsnames]{xcolor}
\usepackage{tikz}
\usetikzlibrary{backgrounds}
\usetikzlibrary{arrows,shapes}
\usepackage{amsmath}
\usepackage{amsthm}
\usepackage{amssymb}
\usepackage{algpseudocode}
\usepackage{mathtools, nccmath}
\usepackage{wrapfig}
\usepackage{comment}
\usepackage{algpseudocode}
\usepackage{makecell}
\usepackage{multirow}
\usepackage[colorlinks=true, linkcolor=red, urlcolor=red]{hyperref}
\usepackage{blindtext}

\usepackage{graphicx}
\usepackage{xspace}

\usepackage{array}
\usepackage{ragged2e}
\newcolumntype{P}[1]{>{\RaggedRight\hspace{0pt}}p{#1}}
\newcolumntype{X}[1]{>{\RaggedRight\hspace*{0pt}}p{#1}}

\usepackage{tcolorbox}

\usepackage{tikz}
\usetikzlibrary{arrows,shapes,positioning,shadows,trees,mindmap}
\usepackage[edges]{forest}
\usetikzlibrary{arrows.meta}
\colorlet{linecol}{black!75}
\usepackage{xkcdcolors} % xkcd colors

\usepackage{tikz}
\usetikzlibrary{backgrounds}
\usetikzlibrary{arrows,shapes}
\usetikzlibrary{tikzmark}
\usetikzlibrary{calc}
\newcommand{\thicktoprule}{\Xhline{4\arrayrulewidth}}
\newcommand{\thickbottomrule}{\Xhline{4\arrayrulewidth}}

\usepackage{lineno,amsmath,amsfonts,amssymb,amsthm,mathrsfs}
\usepackage{hyperref}
\hypersetup{colorlinks=true,linkcolor=blue}
\usepackage[linesnumbered,ruled,vlined,onelanguage]{algorithm2e}
\usepackage{booktabs}
\usepackage[skip=0pt]{caption,subcaption}
\usepackage{float}
\usepackage{longtable}
\usepackage{setspace}
\usepackage{tabularx}
\usepackage{booktabs}
\usepackage{tablefootnote}

\usepackage{todonotes}
\journal{Transportation Research Part B: Methodological}
\biboptions{authoryear}
\newtheorem{theorem}{Theorem}
\newtheorem{lemma}{Lemma}

\newtheorem{definition}{Definition}
\newtheorem{remark}{Remark}
\begin{document}

\begin{frontmatter}

\title{A Scalable Min-Max Multi-Gradient Descent Method for Multi-Objective Transportation Problems}

\author[1]{Yuan-Zheng Lei}
\ead{yzlei@umd.edu}
\author[1]{Yaobang Gong\corref{cor1}}
\ead{ybgong@umd.edu}
\cortext[cor1]{Corresponding author.}
\author[1]{Xianfeng Terry Yang}
\ead{xtyang@umd.edu}
%\author{Authors’ names are not included for peer review}

%\cortext[cor1]{Corresponding author. Xianfeng Terry Yang}
\address[1]{Department of Civil \& Environmental Engineering, University of Maryland, 1173 Glenn Martin Hall, College Park, MD 20742, United States}
\begin{abstract}
\par This paper develops a scalable min-max multi-gradient descent framework for multi-objective transportation problems, motivated by the growing need to optimize complex transportation systems under competing objectives and high-dimensional constraints. Conventional linear scalarization can obscure genuine trade-offs behind pre-specified weights, while evolutionary heuristics often become inefficient when the feasible region is shaped by many constraints. In contrast, the proposed framework directly exploits gradient information to construct a balanced common descent direction. At each iteration, it solves a constrained min-max direction-finding problem that minimizes the worst directional derivative among all objectives, thereby explicitly improving the least-improved objective and promoting robust progress across competing losses or performance criteria. When the problem-specific constraints are linear, this direction-finding problem reduces to a linear program. Theoretically, we prove the existence of the proposed direction, establish its connection to Pareto criticality, and analyze the convergence behavior for a fixed-step variant of the update. We further show that, under a dense worst-case interior-point model, the proposed linear-programming-based formulation has the same asymptotic complexity order as classical quadratic-programming-based steepest-descent formulations, while offering a different and more balanced direction-selection mechanism. Numerical experiments on two representative transportation applications demonstrate the practical value of this framework. In training a physics-informed machine learning-based car-following model, the proposed method eliminates the need for manually fixed loss weights and achieves the best average predictive accuracy among all compared training strategies, reducing the average RMSE by approximately \(5.81\%\) relative to the baseline while substantially lowering variability across random seeds. In a multi-criteria traffic assignment problem balancing user equilibrium and system optimum, the method generates stable common descent directions and delivers favorable performance on large transportation networks. Overall, the proposed framework offers a robust, scalable, and theoretically grounded gradient-based alternative for solving multi-objective transportation problems.
\end{abstract}
\end{frontmatter}

\section{Introduction} \label{1}
%\linenumbers
\par Multi-objective optimization (MO) has long attracted attention due to its fundamental theoretical significance and wide-ranging practical applications. In real-world scenarios, decision-makers are rarely guided by a single criterion; rather, they must balance multiple, often conflicting, objectives (\cite{miettinen1999nonlinear, ehrgott2005multicriteria}). This is particularly evident in the field of transportation. Examples include multi-objective static user equilibrium (\cite{leurent1996theory}, \cite{nagurney2000multiclass}, \cite{nagurney2002multiclass}, and \cite{yang2004multi}) and stochastic user equilibrium problems (\cite{ehrgott2015multi}, \cite{sun2019multi}, \cite{jaber2009mixed}, and \cite{wang2014bi}), multi-criteria shortest path problems (\cite{chen2013bicriterion} and \cite{jafari2017multicriteria}), multi-objective tolled facilities pricing problems (\cite{he2017optimal}), multi-objective routing problems (\cite{pacheco2013bi}, \cite{jozefowiez2008multi}, \cite{kerbache2000multi}, and \cite{coelho2017multi}), and multi-criteria calibration problems (\cite{punzo2021calibration}). Despite the variety of problem types, the modeling and solution of MO problems in transportation research still exhibit certain limitations.

\par First, it is common to address multi-objective problems by applying linear scalarization, which transforms them into single-objective problems through a weighted sum of the original objectives. This approach implicitly assumes a linear trade-off among the objectives, which often fails to reflect real-world preferences or trade-offs. Furthermore, the choice of scalarization weights is typically heuristic and lacks theoretical or empirical justification, a concern that is frequently overlooked in existing studies. Second, evolutionary algorithms (EA) remain the dominant class of solution methods in the transportation literature for multi-objective problems. While these algorithms offer flexibility and require no gradient information, they are not well-suited for problems with a large number of constraints. Even when the constraints are all linear, a large number of constraints can significantly narrow the feasible region, thereby limiting the effectiveness of the population-based stochastic search process. Equality constraints exacerbate this issue by restricting feasible solutions to lower-dimensional manifolds, making random sampling within the feasible region increasingly unlikely. Moreover, due to the absence of gradient-based search directions, evolutionary algorithms scale poorly with problem size and are often confined to small or synthetic instances. As a result, many EA-based studies in transportation deal only with toy problems and cannot be easily extended to realistic, large-scale applications.

\par To address these limitations, this paper proposes a min--max programming approach based on multi-gradient descent algorithms (MGDA), specifically designed for multi-objective optimization problems encountered in transportation systems. The core idea of MGDA is to find a common update direction that improves all objectives to first order whenever such a direction exists. Two most famous and classical MGDA formulations are particularly relevant. \cite{fliege2000steepest} proposed a multi-gradient steepest descent method, which generalizes the single-objective steepest descent algorithm. The common descent direction $\mathbf{d}$ can be obtained by solving the following quadratic programming problem:
\begin{equation}
    \min_{\mathbf{d},\, \gamma} 
    \left\{ 
    \gamma + \frac{1}{2} \|\mathbf{d}\|^2 
    \;\middle|\; 
    \nabla \mathcal{F}_t(\boldsymbol{\theta})^{\top} \mathbf{d} \leq \gamma, 
    \; \forall t = 1, \dots, T 
    \right\}
    \label{eq:1}
\end{equation}
where $\mathcal{F}_t(\cdot)$ are the objective functions and $\boldsymbol{\theta} \in \mathbb{R}^d$ denotes the shared decision variables.

In particular, the minimum norm point method introduced by \cite{desideri2012multiple}, where the minimum norm point of the convex hull of the objective gradients serves as the common descent direction, has attracted considerable attention. Mathematically, the corresponding convex coefficients can be obtained by solving the following linearly constrained quadratic programming problem:
\begin{equation}
    \underset{\mu^{1}, \ldots, \mu^{T}}{\min}
    \left\{
    \left\|
    \sum_{t=1}^{T}\mu^{t}\nabla{\mathcal{F}}_{t}(\boldsymbol{\theta})
    \right\|^{2}
    \;\middle|\;
    \sum_{t=1}^{T}\mu^{t} = 1,\ 
    \mu^{t} \geq 0,\ \forall t
    \right\}
    \label{eq:2}
\end{equation}
after which the minimum norm point can be expressed as
\[
\mathbf{w} = \sum_{t=1}^{T}\mu^{t}_{*}\nabla{\mathcal{F}}_{t}(\boldsymbol{\theta})
\]
and $-\mathbf{w}$ gives a common descent direction when $\mathbf{w}\neq \mathbf 0$. Essentially, \textcolor{blue}{\textbf{Problem}}~\eqref{eq:2} can be interpreted as a dual or equivalent minimum-norm reformulation of the classical steepest descent subproblem in \textcolor{blue}{\textbf{Problem}}~\eqref{eq:1} under the standard unconstrained MGDA setting\footnote{Refer to \textcolor{blue}{\textbf{Theorem}} \ref{theorem:3}.}.

\par By effectively leveraging gradient information, MGDA is generally considered more efficient than evolutionary algorithms when the objective functions are differentiable, smooth, and Lipschitz continuous, particularly in high-dimensional problems. Moreover, MGDA has demonstrated the capability to handle optimization tasks involving millions of variables, as commonly encountered in neural network training (\cite{sener2018multi}). These attributes have made MGDA central in multi-task learning: viewing multi-task learning as a multi-objective optimization problem, \cite{sener2018multi} proposed an MGDA approach based on the Frank--Wolfe algorithm. To address high-dimensional challenges, they introduced an upper bound on the multi-objective loss and optimized this bound under realistic assumptions. Subsequent research has sought more evenly distributed and continuous Pareto fronts. For instance, \cite{lin2019pareto} reformulated the multi-objective problem as several constrained subproblems that capture distinct trade-off preferences; \cite{ma2020efficient} later expanded on this concept to construct locally continuous Pareto sets in contrast to previous methods yielding finite, sparse, or discrete solutions. Likewise, \cite{liu2021profiling} proposed a gradient-based algorithm using Stein variational gradient descent to efficiently discover diverse Pareto solutions in high-dimensional spaces without requiring predefined preferences.

\par Building on these ideas, \cite{phan2022stochastic} addressed probabilistic multi-objective optimization using Stochastic Multiple Target Sampling Gradient Descent (MT-SGD). Further advancing this line of research, \cite{momma2022multi} proposed a generic framework for identifying Pareto-optimal solutions across various preference structures, and \cite{mercier2018stochastic} extended both stochastic gradient descent (SGD) and MGDA to objective functions expressed as expectations of random functions. Numerous other studies have similarly demonstrated MGDA's versatility in domains such as recommendation systems (\cite{lin2019apareto, milojkovic2019multi, chen2021pareto, pan2023fedmdfg}). Recent work adopting MGDA emphasizes generating comprehensive, continuous Pareto fronts or those aligning with particular user preferences (\cite{mahapatra2020multi}).

\par Despite its wide applicability, most well-known MGDA frameworks (\cite{fliege2000steepest}, \cite{schaffler2002stochastic}, and \cite{desideri2012multiple}) focus on unconstrained multi-objective optimization. In constrained scenarios, the minimum norm point derived from the convex hull of only the objective gradients can differ from the direction relevant to the feasible set, which may forfeit the mathematical guarantee that the resulting direction remains a valid feasible common descent direction. For this reason, the approach mentioned in \cite{desideri2012multiple} cannot be applied directly to CMO problems. Other MGDA algorithms, such as \cite{fliege2000steepest}, \cite{fliege2009newton}, and \cite{ansary2021globally}, can handle constraints or be adapted to them. For example, \cite{ansary2021globally} developed a variant in which the common descent direction is obtained by solving a quadratically constrained quadratic programming problem, a natural extension of \cite{fliege2000steepest}; \cite{fliege2009newton} proposed a Newton-type method for the unconstrained case, where the common descent direction is obtained by solving a quadratic programming problem (QP). In short, most existing MGDA methods that are usable for CMO ultimately rely on QP subproblems. Nevertheless, a robust and balanced constrained MGDA framework remains a valuable direction of investigation.

\par The classical steepest descent variant of \cite{fliege2000steepest} can be interpreted as a regularized min-max problem, in which the objective minimizes an upper bound on all directional derivatives while simultaneously penalizing the step norm. This norm regularization promotes smoother and more stable updates, but it also changes the direction-selection criterion by balancing worst-case descent against the norm penalty. In contrast, a pure min--max formulation directly minimizes the maximum directional derivative, thereby focusing explicitly on improving the least-improved objective at each iteration. This property is particularly attractive when the goal is to obtain a balanced common descent direction across competing objectives or loss terms.

\par As noted above, most MGDA approaches usable in MO rely on QP. By contrast, the proposed min--max formulation can be equivalently transformed into a linear program (LP). In constrained multi-objective optimization (CMO) problems, problem-specific constraints are present; when these constraints are also linear, the search for a common descent direction reduces to a purely linear programming problem, in contrast to the QP used traditionally. While LP is structurally simpler than QP, the formulation may introduce additional constraints, so standard interior-point complexity analysis does not confer an inherent worst-case advantage to LP. In other words, while the LP-based formulation does not necessarily improve the per-iteration worst-case computational complexity, it provides a different and more balanced common descent direction, which can lead to better practical convergence and improved final solutions.

\par Accordingly, we propose a min-max MGDA framework for handling multi-objective transportation problems. The method focuses on optimizing the worst-case performance of the common descent direction while ensuring that the direction remains valid for all objectives and respects the problem-specific constraints. This construction emphasizes robustness and balance in direction selection under the common descent requirement. The primary contributions of our work are as follows:
\begin{itemize}
\item We propose a min-max programming algorithm designed to determine a robust and balanced common descent direction under convex constraints. The algorithm theoretically guarantees Pareto criticality. Distinct from existing multi-gradient descent methods relying on quadratic programming, the proposed direction-finding subproblem is formulated as a linear program when problem-specific constraints are linear, and directly optimizing the worst directional derivative to achieve a balanced trade-off among objectives.

\item We demonstrate the effectiveness of the proposed framework in training a physics-informed machine-learning-based car-following model. In this application, the data-fitting loss and the physics-based loss are treated as separate objectives, and the proposed method computes a common descent direction without manually fixing scalarization weights. The numerical results show that our method achieves the best final predictive performance among all compared training strategies, reducing the average RMSE by approximately \(5.81\%\) relative to the baseline while also reducing the variability across random seeds.

\item We further evaluate the proposed method on a large-scale multi-criteria traffic assignment problem that balances user equilibrium (UE) and system optimum (SO) objectives. The results show that the proposed min-max programming method exhibits better convergence behavior than the classical steepest descent formulation: it reaches high-quality solutions with substantially fewer iterations and lower total CPU time on the tested transportation networks.

\item We provide theoretical analysis and empirical validation of the computational complexity of our linear-programming-based formulation. Although it does not yield an inherent dense worst-case complexity advantage over QP-based steepest-descent formulations, it produces robust common descent directions and achieves favorable practical performance on large-scale problems with sparse linear constraints and without any constraints.
\end{itemize}

\section{Notation}
\par This section lists the notation and symbols consistently used in the remainder of the paper, as summarized in \textcolor{blue}{\textbf{Table}} \ref{table:1}.
\begin{table}[htbp]
  \caption{\label{table:1}Notations and symbols}
  \centering
  \begin{tabularx}{\textwidth}{@{} l X @{}}
    \thicktoprule
    \textbf{Symbols} & \textbf{Descriptions} \\
    \midrule
    $\|\cdot\|$ &  
    The Euclidean ($\ell_2$) norm, i.e., $\|x\|_2 = \big(\sum_i x_i^2\big)^{1/2}$\\
    $\langle x, y \rangle$ & 
    The inner product of $x$ and $y$; in $\mathbb{R}^n$, $\langle x, y \rangle = \sum_{i=1}^n x_i y_i$ \\
    $f_i(\cdot)$ & Problem-specific constraints, i.e., constraints in CMO problems, index by $i, i = 1,..,m$ \\
    $\mathcal{F}_{t}(\cdot)$ & Objective functions, indexed by $t,t = 1,...,T$ \\
    $\boldsymbol{\theta}$ &     $\boldsymbol{\theta}$ is the vectorized decision variable set \\
    $\asymp(\cdot)$ & $\asymp(\cdot)$ is the asymptotic complexity notation. 
    $f(n) = \asymp(g(n))$ means that $f(n)$ grows at the same order as $g(n)$, 
    i.e., there exist constants $c_1,c_2>0$ and $n_0$ such that 
    $c_1 g(n) \le f(n) \le c_2 g(n)$ for all $n \ge n_0$ \\
    $\mathbf{\Theta}$   & $\mathbf{\Theta}\subseteq\mathbb{R}^{n}$ is the solution domain\\
    $\eta$ & $\eta$ is an auxiliary variable introduced to reformulate the original $\min\max$ optimization problem as a linear program, by constraining all objective functions to be bounded under $\eta$ \\
    $\beta$ & $\beta$ is the decay factor in the step size search algorithm, with $\beta\in(0,1)$ \\
    $\mathbf{d}$ & $\mathbf{d}$ is both a direction vector and a variable; $\mathbf{d}\in\mathbb{R}^{\dim(\boldsymbol{\theta})}$ and matches the structure of $\boldsymbol{\theta}$ \\
    $h$ & $h$ refers to the step size \\
    $\phi$ & $\phi$ is a small positive number that ensures the feasibility of the linear approximation \\
    $\triangleq$ & denotes a definition, i.e., is defined as \\
    $\downarrow$ & $\alpha \downarrow 0$ means that $\alpha$ approaches $0$ from the positive side, i.e., $\alpha \to 0^{+}$ \\
    $o(\cdot)$ & For two functions $r(\alpha)$ and $\alpha$, the relation $r(\alpha)=o(\alpha)$ as $\alpha\to 0^{+}$ means that $\displaystyle \lim_{\alpha\to 0^{+}} \frac{r(\alpha)}{\alpha}=0$ \\
    \thickbottomrule
  \end{tabularx}
\end{table}
\section{Scalable min-max programming for unconstrained/constrained multi-gradient descent}
\par Before formally describing the min-max approach, we will review important concepts and outline some key assumptions. We consider the following general constrained multi-objective optimization problem:
\begin{subequations}\label{eq:3}
\begin{align}
\min_{\boldsymbol{\theta}}
  &\;\bigl(\mathcal{F}_{1}(\boldsymbol{\theta}),\dots,\mathcal{F}_{T}(\boldsymbol{\theta})\bigr)
  \label{eq:3a}\\
\text{s.t.}\quad
  &\;f_i(\boldsymbol{\theta}) \le 0,\quad i=1,2,\dots,m
  \label{eq:3b}
\end{align}
\end{subequations}
where $\mathcal{F}(\cdot)$ denotes the objective functions and $f(\cdot)$ denotes the constraint functions.  The following assumptions are made:
\begin{itemize}
\item Each \(\mathcal{F}_{t} : \mathbb{R}^{n} \to \mathbb{R}\) is differentiable everywhere on its domain \(\mathbf{\Theta}\), for all \(t = 1, \dots, T\).
\item  \(f_{i} : \mathbb{R}^{n} \to \mathbb{R}\) are convex, for $i = 1,...,m$.
\item  There exists a point \(x_{s} \in \mathbb{R}^{n}\)  such that
\[
    f_i(x_{s}) \;<\; 0 
    \quad \text{for all } i=1,\dots,m.
\]
[Slater’s Condition (\cite{ben2004lecture})]
\end{itemize}
\begin{definition}[Pareto dominance] \label{definition:1}
Under the same problem setting, a solution $\boldsymbol{\theta}$ dominates a solution $\hat{\boldsymbol{\theta}}$ if and only if $\mathcal{F}_{t}(\boldsymbol{\theta}) \leq \mathcal{F}_{t}(\hat{\boldsymbol{\theta}}), \forall t \in   \{t = 1, \dots, T\}$, and $\mathcal{F}_{t'}(\boldsymbol{\theta}) < \mathcal{F}_{t'}(\hat{\boldsymbol{\theta}}), \exists t' \in  \{t = 1, \dots, T\}$.
\end{definition}
\begin{definition}[Pareto optimal] \label{definition:2} A solution $\boldsymbol{\theta}^{*}$ is called Pareto optimal if and only if there exists no solution $\boldsymbol{\theta}$ that dominates $\boldsymbol{\theta}^{*}$.
\end{definition}
\begin{definition}[Pareto critical point (\cite{fliege2000steepest})]
\label{definition:3}
Under the same problem setting, a feasible point $\boldsymbol{\theta}^{*}$ is called a Pareto critical point if there is no nonzero feasible direction $\mathbf d$ such that
\begin{equation}
\nabla \mathcal F_t(\boldsymbol{\theta}^{*})^{\top}\mathbf d < 0
\qquad t=1,\ldots,T
\label{eq:4}
\end{equation}
and
\begin{equation}
\nabla f_i(\boldsymbol{\theta}^{*})^{\top}\mathbf d \le 0
\qquad \forall i\in 1,...,m
\label{eq:5}
\end{equation}
Equivalently, $\boldsymbol{\theta}^{*}$ is a Pareto critical point\footnote{In our definition, we essentially consider a constrained version of the Pareto critical point compared to \cite{fliege2000steepest} where extra convex constraints $f_{i}(\cdot),i = 1,..,m$ are introduced.} if and only if every feasible direction $\mathbf{d} \neq \mathbf 0$ fails to decrease all objectives simultaneously to first order.
\end{definition}
\par We propose a min-max search algorithm to identify a common descent direction that ultimately achieves Pareto critical.  
\par Essentially, we seek a common descent direction by solving the following optimization problem:
\begin{subequations}\label{eq:6}
\begin{align}
  &w
    = \min_{\mathbf{d}}\;\max_{t=1,2,\dots,T}
      \,\nabla\mathcal{F}_{t}(\boldsymbol{\theta})^{\top}\mathbf{d}
    \label{eq:6a}\\
  &\text{s.t.}\quad
    f_{i}(\boldsymbol{\theta} + \mathbf{d}) \le 0,\quad i=1,2,\dots,m
    \label{eq:6b}\\
  &\quad
    \nabla \mathcal{F}_{t}(\boldsymbol{\theta})^{\top}\mathbf{d} \le 0,
    \quad t=1,2,\dots,T
    \label{eq:6c}\\
  &\quad
    -\mathbf{1} \;\le\; \mathbf{d} \;\le\; \mathbf{1}
    \label{eq:6d}
\end{align}
\end{subequations}
\par In \textcolor{blue}{\textbf{Problem}} \ref{eq:6}, the goal is to identify a direction that minimizes the upper bound of all directional derivatives while satisfying the given constraints. The first set of constraints ensure that the candidate direction \(\textbf{d}\) remains feasible for a step size of at least one, leaving ample flexibility for choosing step sizes later on, since the condition \( \nabla \mathcal{F}_{t}(\boldsymbol{\theta})^{\top} \mathbf{d} \leq 0 \) ensures a common descent direction, we adopt a step size searching approach as described by \cite{desideri2012multiple}. Specifically, we aim to find the maximum monotonic decreasing interval \( [\boldsymbol{\theta}, \boldsymbol{\theta} + h \cdot \mathbf{d}^{*}] \) for all objectives, where \( h \) will serve as the final update step length, mathematically:
 \begin{equation}
h = \max \{ h' > 0 \mid \nabla \mathcal{F}_{t}(\boldsymbol{\theta} + \tau \mathbf{d}^{*})^{\top} \mathbf{d}^{*} \leq 0, \, \forall \tau \in (0, h'), \, t = 1, 2, \ldots, T \} \label{eq:7}
\end{equation}
where $\mathbf{d}^{*}$ is the common descent direction obtained by solving \textcolor{blue}{\textbf{Problem}} \ref{eq:6}. The process of finding $h$ is a line search process, during which the value of $h$ may constantly decrease until a maximum monotonic decreasing interval is found.
\par In our experiments, we discovered that maintaining a sufficiently large feasible step size is crucial for the step searching process. Specifically, we impose a stringent constraint \ref{eq:6b}, which ensures that any chosen direction $\mathbf{d}$ remains feasible at a full step length. 
\par Without this explicit constraint, a selected direction may satisfy feasibility only for an infinitesimally small step size. Mathematically, suppose we select a direction $\mathbf{d}$ without directly enforcing the step size constraint. Then, the largest feasible step size along $\mathbf{d}$, defined as
\begin{equation} \label{eq:8}
h_{\max} = \max\{h>0\,|\,f_i(\boldsymbol{\theta}+h\mathbf{d})\leq 0,\,i=1,\dots,m\}
\end{equation}
can become arbitrarily small, effectively limiting the practical progress of parameter updates. By explicitly enforcing feasibility at a relatively large step size, we sacrifice some freedom in the feasible region to guarantee that every selected direction provides meaningful descent at a substantial, predetermined step size, thus significantly enhancing optimization efficiency and stability. 
\par In addition, the condition $\nabla \mathcal{F}_{t}(\boldsymbol{\theta})^{\top} \mathbf{d} \leq 0, \quad t = 1,2,...,t$ forces \(\textbf{d}\) to be a common descent direction, and $\mathbf{-1} \leq \mathbf{d}  \leq \mathbf{1}$ ensures the search direction is bounded. 
\par Before presenting the convergence analysis, we first establish the Pareto-criticality interpretation of \textcolor{blue}{\textbf{Problem}}~\eqref{eq:6}. The quantity \(w(\boldsymbol{\theta})\) measures whether there exists a feasible direction that strictly decreases all objectives to first order. If \(w(\boldsymbol{\theta})<0\), the current point admits a feasible common descent direction. If \(w(\boldsymbol{\theta})=0\), no such feasible common descent direction exists, and the current point is Pareto critical. Therefore, \(w(\boldsymbol{\theta})\) provides both a direction-generation criterion and a natural stopping criterion for the proposed min--max programming framework.
\begin{theorem}\label{theorem:1}
Assume that each objective function $\mathcal F_t$ is continuously differentiable, each constraint function $f_i$ is convex and continuously differentiable, and Slater's condition holds, namely, there exists $\boldsymbol{\theta}^{\mathrm s}\in\mathbb R^n$ such that
\begin{equation}
f_i(\boldsymbol{\theta}^{\mathrm s})<0
\qquad i=1,\ldots,m
\label{eq:9}
\end{equation}
For a feasible point $\boldsymbol{\theta}$, let
\begin{equation}
w(\boldsymbol{\theta})
\triangleq
\min_{\mathbf d}
\max_{t=1,\ldots,T}
\nabla \mathcal F_t(\boldsymbol{\theta})^\top \mathbf d
\label{eq:10}
\end{equation}
subject to the constraints in \textcolor{blue}{\textbf{Problem}}~\eqref{eq:6}. Then the following statements hold:
\begin{itemize}
    \item[(i)] \(w(\boldsymbol{\theta})\le 0\).
    \item[(ii)] If \(w(\boldsymbol{\theta})<0\), then any optimal solution \(\mathbf d^*\) of \textcolor{blue}{\textbf{Problem}}~\ref{eq:6} satisfies
    \begin{equation}
    \nabla \mathcal F_t(\boldsymbol{\theta})^\top \mathbf d^*
    \le
    w(\boldsymbol{\theta})
    <
    0
    \qquad t=1,\ldots,T
    \label{eq:11}
    \end{equation}
    and hence \(\mathbf d^*\) is a common descent direction that strictly decreases all objectives to first order.
    \item[(iii)] If \(w(\boldsymbol{\theta})=0\), then \(\boldsymbol{\theta}\) is a Pareto critical point in the sense of \textcolor{blue}{\textbf{Definition}}~\ref{definition:3}.
\end{itemize}
\end{theorem}

\begin{proof}
Let
\begin{equation}
\mathcal D(\boldsymbol{\theta})
\triangleq
\Bigl\{
\mathbf d\in\mathbb R^n:
\boldsymbol{\theta}+\mathbf d\in\mathcal X,\ 
\nabla \mathcal F_t(\boldsymbol{\theta})^\top \mathbf d\le 0,\ t=1,\ldots,T,\ 
-\mathbf 1\le \mathbf d\le \mathbf 1
\Bigr\}
\label{eq:12}
\end{equation}
Since \(\mathbf d=\mathbf 0\) satisfies
\begin{equation}
\boldsymbol{\theta}+\mathbf 0=\boldsymbol{\theta}\in\mathcal X
\qquad
\nabla \mathcal F_t(\boldsymbol{\theta})^\top \mathbf 0=0,\ t=1,\ldots,T
\qquad
-\mathbf 1\le \mathbf 0\le \mathbf 1
\label{eq:13}
\end{equation}
the set \(\mathcal D(\boldsymbol{\theta})\) is nonempty. Moreover, \(\mathcal D(\boldsymbol{\theta})\) is a closed subset of the compact box \([ -1,1]^n\), because \(\mathcal X\) is closed by continuity of the constraint functions, the inequalities
\begin{equation}
\nabla \mathcal F_t(\boldsymbol{\theta})^\top \mathbf d\le 0
\qquad t=1,\ldots,T
\label{eq:14}
\end{equation}
define closed half-spaces, and the box constraints are closed. Hence \(\mathcal D(\boldsymbol{\theta})\) is compact. Since the map
\begin{equation}
\mathbf{d}  \mapsto \max_{t=1,\ldots,T}\nabla \mathcal F_t(\boldsymbol{\theta})^\top \mathbf d
\label{eq:15}
\end{equation}
is continuous\footnote{For fixed $\boldsymbol{\theta}$, each mapping $\mathbf d \mapsto \nabla \mathcal F_t(\boldsymbol{\theta})^\top \mathbf d$ is linear in $\mathbf d$, hence continuous. Since the pointwise maximum of finitely many continuous functions is continuous, the map $\mathbf d \mapsto \max_{t=1,\ldots,T}\nabla \mathcal F_t(\boldsymbol{\theta})^\top \mathbf d$ is continuous.}, the minimum in \textcolor{blue}{\textbf{Problem}} \eqref{eq:6} is attained on \(\mathcal D(\boldsymbol{\theta})\). Therefore an optimal solution \(\mathbf d^*\) exists.

Because \(\mathbf d=\mathbf 0\) is feasible for \textcolor{blue}{\textbf{Problem}}~\eqref{eq:6}, its objective value is
\begin{equation}
\max_{t=1,\ldots,T}
\nabla \mathcal F_t(\boldsymbol{\theta})^\top \mathbf 0
=
0
\label{eq:16}
\end{equation}
Therefore, the optimal value, being the minimum of the objective over all feasible directions, satisfies
\begin{equation}
w(\boldsymbol{\theta})\le 0
\label{eq:17}
\end{equation}
This proves (i).

Let \(\mathbf d^*\) be any optimal solution of \textcolor{blue}{\textbf{Problem}}~\eqref{eq:6}. By the definition of the optimal value, we have
\begin{equation}
\max_{t=1,\ldots,T}
\nabla \mathcal F_t(\boldsymbol{\theta})^\top \mathbf d^*
=
w(\boldsymbol{\theta})
\label{eq:18}
\end{equation}
If \(w(\boldsymbol{\theta})<0\), then \eqref{eq:18} immediately implies
\begin{equation}
\nabla \mathcal F_t(\boldsymbol{\theta})^\top \mathbf d^*
\le
w(\boldsymbol{\theta})
<
0
\qquad t=1,\ldots,T
\label{eq:19}
\end{equation}
Thus \(\mathbf d^*\) strictly decreases every objective to first order, proving (ii).

Now assume
\begin{equation}
w(\boldsymbol{\theta})=0
\label{eq:20}
\end{equation}
Suppose, for contradiction, that \(\boldsymbol{\theta}\) is not Pareto critical in the sense of \textcolor{blue}{\textbf{Definition}}~\ref{definition:3}. Then there exists a nonzero feasible direction \(\mathbf p\in\mathbb R^n\) such that
\begin{equation}
\nabla \mathcal F_t(\boldsymbol{\theta})^\top \mathbf p<0
\qquad t=1,\ldots,T
\label{eq:21}
\end{equation}
and
\begin{equation}
\nabla f_i(\boldsymbol{\theta})^\top \mathbf p\le 0
\qquad \forall i\in I(\boldsymbol{\theta})
\label{eq:22}
\end{equation}
where
\begin{equation}
I(\boldsymbol{\theta})
\triangleq
\{\, i\in\{1,\ldots,m\}: f_i(\boldsymbol{\theta})=0 \,\}
\label{eq:23}
\end{equation}

Let
\begin{equation}
\mathbf v
\triangleq
\boldsymbol{\theta}^{\mathrm s}-\boldsymbol{\theta}
\label{eq:24}
\end{equation}
where \(\boldsymbol{\theta}^{\mathrm s}\) is the Slater point in \eqref{eq:9}. For any \(i\in I(\boldsymbol{\theta})\)\footnote{$I(\boldsymbol{\theta})$ is the active set.}, convexity and differentiability of \(f_i\) give
\begin{equation}
f_i(\boldsymbol{\theta}^{\mathrm s})
\ge
f_i(\boldsymbol{\theta})
+
\nabla f_i(\boldsymbol{\theta})^\top(\boldsymbol{\theta}^{\mathrm s}-\boldsymbol{\theta})
=
\nabla f_i(\boldsymbol{\theta})^\top \mathbf v
\label{eq:25}
\end{equation}
because \(f_i(\boldsymbol{\theta})=0\) for all \(i\in I(\boldsymbol{\theta})\). Since \(f_i(\boldsymbol{\theta}^{\mathrm s})<0\), it follows that
\begin{equation}
\nabla f_i(\boldsymbol{\theta})^\top \mathbf v<0
\qquad \forall i\in I(\boldsymbol{\theta})
\label{eq:26}
\end{equation}

For \(\lambda>0\), define
\begin{equation}
\mathbf q_\lambda
\triangleq
\mathbf p+\lambda \mathbf v
\label{eq:27}
\end{equation}
Then, by \eqref{eq:22} and \eqref{eq:26},
\begin{equation}
\nabla f_i(\boldsymbol{\theta})^\top \mathbf q_\lambda
=
\nabla f_i(\boldsymbol{\theta})^\top \mathbf p
+
\lambda \nabla f_i(\boldsymbol{\theta})^\top \mathbf v
<
0
\qquad \forall i\in I(\boldsymbol{\theta})
\label{eq:28}
\end{equation}
Moreover, since \(\nabla \mathcal F_t(\boldsymbol{\theta})^\top \mathbf p<0\) for all \(t\), continuity with respect to \(\lambda\) implies that there exists \(\bar\lambda>0\) such that for every \(\lambda\in(0,\bar\lambda]\),
\begin{equation}
\nabla \mathcal F_t(\boldsymbol{\theta})^\top \mathbf q_\lambda<0
\qquad t=1,\ldots,T
\label{eq:29}
\end{equation}
Fix such a \(\lambda\), and write
\begin{equation}
\mathbf q\triangleq \mathbf q_\lambda
\label{eq:30}
\end{equation}

For any active constraint \(i\in I(\boldsymbol{\theta})\), differentiability of \(f_i\) gives
\begin{equation}
f_i(\boldsymbol{\theta}+\alpha\mathbf q)
=
f_i(\boldsymbol{\theta})
+
\alpha \nabla f_i(\boldsymbol{\theta})^\top \mathbf q
+
o(\alpha)
=
\alpha \nabla f_i(\boldsymbol{\theta})^\top \mathbf q
+
o(\alpha)
\qquad \alpha\downarrow 0
\label{eq:31}
\end{equation}
Since \(\nabla f_i(\boldsymbol{\theta})^\top \mathbf q<0\) by \eqref{eq:28}, there exists \(\alpha_i>0\) such that
\begin{equation}
f_i(\boldsymbol{\theta}+\alpha\mathbf q)<0
\qquad \forall \alpha\in(0,\alpha_i]
\label{eq:32}
\end{equation}
For any inactive constraint \(i\notin I(\boldsymbol{\theta})\), we already have \(f_i(\boldsymbol{\theta})<0\), so continuity yields the existence of \(\beta_i>0\) such that
\begin{equation}
f_i(\boldsymbol{\theta}+\alpha\mathbf q)<0
\qquad \forall \alpha\in(0,\beta_i]
\label{eq:33}
\end{equation}
With the convention that the minimum over an empty index set is \(+\infty\), define
\begin{equation}
\bar\alpha
\triangleq
\min\Bigl\{
1,\ 
\min_{i\in I(\boldsymbol{\theta})}\alpha_i,\ 
\min_{i\notin I(\boldsymbol{\theta})}\beta_i,\ 
\frac{1}{\|\mathbf q\|_\infty}
\Bigr\}
\label{eq:34}
\end{equation}
Choose any \(\alpha\in(0,\bar\alpha]\), and define
\begin{equation}
\mathbf d\triangleq \alpha \mathbf q
\label{eq:35}
\end{equation}
Then, by construction,
\begin{equation}
\boldsymbol{\theta}+\mathbf d\in\mathcal X
\label{eq:36}
\end{equation}
and
\begin{equation}
-\mathbf 1\le \mathbf d\le \mathbf 1
\label{eq:37}
\end{equation}
Furthermore, by \eqref{eq:29},
\begin{equation}
\nabla \mathcal F_t(\boldsymbol{\theta})^\top \mathbf d
=
\alpha \nabla \mathcal F_t(\boldsymbol{\theta})^\top \mathbf q
<
0
\qquad t=1,\ldots,T
\label{eq:38}
\end{equation}
Thus \(\mathbf d\) is feasible for \textcolor{blue}{\textbf{Problem}}~\eqref{eq:6}, and its objective value satisfies
\begin{equation}
\max_{t=1,\ldots,T}
\nabla \mathcal F_t(\boldsymbol{\theta})^\top \mathbf d
<
0
\label{eq:39}
\end{equation}
Hence \(w(\boldsymbol{\theta})<0\), contradicting \eqref{eq:20}. Therefore \(\boldsymbol{\theta}\) must be Pareto critical. This proves (iii).
\end{proof}

\begin{lemma}\label{lemma:1}
Under the assumptions of \textcolor{blue}{\textbf{Theorem}}~\ref{theorem:1}, a feasible point \(\boldsymbol{\theta}\) is Pareto critical in the sense of \textcolor{blue}{\textbf{Definition}}~\ref{definition:3} if and only if
\begin{equation}
w(\boldsymbol{\theta})=0
\label{eq:40}
\end{equation}
\end{lemma}

\begin{proof}
The implication
\[
w(\boldsymbol{\theta})=0
\quad\Longrightarrow\quad
\boldsymbol{\theta}\ \text{is Pareto critical}
\]
follows directly from Theorem~\ref{theorem:1}(iii).

Conversely, assume that \(\boldsymbol{\theta}\) is Pareto critical. Since \(\mathbf 0\) is feasible for \textcolor{blue}{\textbf{Problem}}~\ref{eq:6}, Theorem~\ref{theorem:1}(i) gives
\begin{equation}
w(\boldsymbol{\theta})\le 0.
\label{eq:41}
\end{equation}
Suppose, for contradiction, that
\begin{equation}
w(\boldsymbol{\theta})<0.
\label{eq:42}
\end{equation}
Then, by Theorem~\ref{theorem:1}(ii), there exists an optimal solution \(\mathbf d^*\) such that
\begin{equation}
\nabla \mathcal F_t(\boldsymbol{\theta})^\top \mathbf d^*<0
\qquad t=1,\ldots,T
\label{eq:43}
\end{equation}
Moreover, \(\mathbf d^*\neq \mathbf 0\), since the inequalities in \eqref{eq:43} cannot hold for \(\mathbf d^*=\mathbf 0\). And since \(\mathbf d^*\) is feasible for \textcolor{blue}{\textbf{Problem}}~\ref{eq:6}, we have
\begin{equation}
\boldsymbol{\theta}+\mathbf d^*\in\mathcal X
\label{eq:44}
\end{equation}
Because \(\mathcal X\) is convex, \eqref{eq:44} implies
\begin{equation}
\boldsymbol{\theta}+\tau \mathbf d^*\in\mathcal X
\qquad \forall \tau\in[0,1]
\label{eq:45}
\end{equation}
and hence \(\mathbf d^*\) is a feasible direction at \(\boldsymbol{\theta}\). Now let \(i\in I(\boldsymbol{\theta})\). By convexity and differentiability of \(f_i\),
\begin{equation}
f_i(\boldsymbol{\theta}+\mathbf d^*)
\ge
f_i(\boldsymbol{\theta})
+
\nabla f_i(\boldsymbol{\theta})^\top \mathbf d^*
=
\nabla f_i(\boldsymbol{\theta})^\top \mathbf d^*
\label{eq:46}
\end{equation}
because \(f_i(\boldsymbol{\theta})=0\) for \(i\in I(\boldsymbol{\theta})\). Since \(f_i(\boldsymbol{\theta}+\mathbf d^*)\le 0\), it follows that
\begin{equation}
\nabla f_i(\boldsymbol{\theta})^\top \mathbf d^*\le 0
\qquad \forall i\in I(\boldsymbol{\theta})
\label{eq:47}
\end{equation}
Combining \eqref{eq:43}, \eqref{eq:45}, and \eqref{eq:47}, we conclude that \(\mathbf d^*\) is a nonzero feasible direction that strictly decreases all objectives, contradicting Pareto criticality. Hence \(w(\boldsymbol{\theta})<0\) is impossible, and together with \eqref{eq:41} this yields
\begin{equation}
w(\boldsymbol{\theta})=0
\label{eq:48}
\end{equation}
This completes the proof.
\end{proof}

\par \textcolor{blue}{\textbf{Theorem}}~\ref{theorem:1} and \textcolor{blue}{\textbf{Lemma}}~\ref{lemma:1} establish the Pareto-criticality interpretation of \textcolor{blue}{\textbf{Problem}} \eqref{eq:6}. We next provide a complementary convergence-rate result for an idealized fixed-step variant of the proposed update. This result should be understood as a standard descent-type analysis for the aggregated objective $\Phi(\boldsymbol{\theta})
=
\sum_{t=1}^{T}\mathcal F_t(\boldsymbol{\theta})
$ rather than as a convergence-rate proof for the exact step-search rule in \textcolor{blue}{\textbf{Algorithm}}~\ref{algorithm:2} or for the adaptive optimizer used in the machine-learning experiment (refer to \textcolor{blue}{\textbf{Section}} \ref{sec:5.1}). Its purpose is to show that, when the generated direction satisfies a sufficient descent condition and remains well aligned with the negative aggregated gradient, the resulting fixed-step sequence ensures a standard sublinear convergence guarantee.

\begin{theorem}
\label{theorem:2}
Define the aggregated objective
\begin{equation}
\Phi(\boldsymbol\theta)
\triangleq
\sum_{t=1}^T \mathcal F_t(\boldsymbol\theta)
\label{eq:49}
\end{equation}
Let $\boldsymbol\theta_0\in\mathcal X$, and let $\{\boldsymbol\theta_k\}$ be generated as follows.

At iteration $k$, let $\mathbf d_k$ be an exact minimizer of \textcolor{blue}{\textbf{Problem}}~\eqref{eq:6} at $\boldsymbol\theta_k$, and update
\begin{equation}
\boldsymbol\theta_{k+1}
=
\boldsymbol\theta_k+h\mathbf d_k
\label{eq:50}
\end{equation}
where the step size $h$ is a fixed constant satisfying
\begin{equation}
0<h\le \min\Bigl\{1,\frac{1}{2L}\Bigr\}
\label{eq:51}
\end{equation}

Assume that

\begin{itemize}
    \item[(A1)] each objective function $\mathcal F_t$ is continuously differentiable on an open neighborhood of $\mathcal X$
    \item[(A2)] each constraint function $f_i$ is convex and continuously differentiable on $\mathbb R^n$
    \item[(A3)] Slater's condition holds, namely, there exists $\boldsymbol\theta^{\mathrm s}\in\mathbb R^n$ such that
    \begin{equation}
    f_i(\boldsymbol\theta^{\mathrm s})<0,
    \qquad i=1,\ldots,m
    \label{eq:52}
    \end{equation}
    \item[(A4)] $\Phi$ is continuously differentiable and $\nabla \Phi$ is $L$-Lipschitz continuous with $L>0$ on an open neighborhood containing all line segments joining $\boldsymbol\theta_k$ and $\boldsymbol\theta_k+h\mathbf d_k$
    \item[(A5)] $\Phi$ is bounded below on $\mathcal X$, and
    \begin{equation}
    \Phi_*
    \triangleq
    \inf_{\boldsymbol\theta\in\mathcal X}\Phi(\boldsymbol\theta)
    >
    -\infty
    \label{eq:53}
    \end{equation}
    \item[(A6)] there exists a constant $M>0$ such that for all $k\ge 0$
    \begin{equation}
    2\langle \nabla \Phi(\boldsymbol\theta_k),\mathbf d_k\rangle
    +
    \|\mathbf d_k\|^2
    \le
    0
    \label{eq:54}
    \end{equation}
    \begin{equation}
    \|\mathbf d_k\|
    \ge
    M\|\nabla \Phi(\boldsymbol\theta_k)\|
    \label{eq:55}
    \end{equation}
\end{itemize}

Then the following statements hold:

\begin{itemize}
    \item[(i)] $\boldsymbol\theta_k\in\mathcal X$ for all $k\ge 0$
    \item[(ii)] if $\mathbf d_k=\mathbf 0$ for some $k$, then $\boldsymbol\theta_k$ is a Pareto critical point in the sense of \textcolor{blue}{\textbf{Definition}}~\ref{definition:3}
    \item[(iii)] for every $K\ge 0$
    \begin{equation}
    \frac{1}{K+1}\sum_{k=0}^K \|\nabla \Phi(\boldsymbol\theta_k)\|^2
    \le
    \frac{4\bigl(\Phi(\boldsymbol\theta_0)-\Phi_*\bigr)}{hM^2(K+1)}
    \label{eq:56}
    \end{equation}
    Consequently,
    \begin{equation}
    \min_{0\le k\le K}\|\nabla \Phi(\boldsymbol\theta_k)\|
    \le
    \sqrt{
    \frac{4\bigl(\Phi(\boldsymbol\theta_0)-\Phi_*\bigr)}{hM^2(K+1)}
    }
    \label{eq:57}
    \end{equation}
\end{itemize}
\end{theorem}

\begin{proof}
We first prove item (i) by induction. Since $\boldsymbol\theta_0\in\mathcal X$ by assumption, the claim holds for $k=0$. Suppose that $\boldsymbol\theta_k\in\mathcal X$. Since $\mathbf d_k$ is feasible for \textcolor{blue}{\textbf{Problem}}~\eqref{eq:6} at $\boldsymbol\theta_k$, we have
\begin{equation}
\boldsymbol\theta_k+\mathbf d_k\in\mathcal X
\label{eq:58}
\end{equation}
Because $0<h\le 1$, the update \eqref{eq:50} can be written as
\begin{equation}
\boldsymbol\theta_{k+1}
=
(1-h)\boldsymbol\theta_k+h(\boldsymbol\theta_k+\mathbf d_k)
\label{eq:59}
\end{equation}
Since $\mathcal X$ is convex by (A2), \eqref{eq:58} and \eqref{eq:59} imply
\begin{equation}
\boldsymbol\theta_{k+1}\in\mathcal X
\label{eq:60}
\end{equation}
Therefore, by induction, $\boldsymbol\theta_k\in\mathcal X$ for all $k\ge 0$. This proves item (i).

Next suppose that $\mathbf d_k=\mathbf 0$ for some $k$. Since $\mathbf d_k$ is an optimal solution of \textcolor{blue}{\textbf{Problem}}~\eqref{eq:6}, its objective value is
\begin{equation}
\max_{t=1,\ldots,T}\nabla \mathcal F_t(\boldsymbol\theta_k)^\top \mathbf d_k
=
0
\label{eq:61}
\end{equation}
Hence, the optimal value satisfies
\begin{equation}
w(\boldsymbol\theta_k)=0
\label{eq:62}
\end{equation}
By \textcolor{blue}{\textbf{Lemma}}~\ref{lemma:1}, \eqref{eq:62} implies that $\boldsymbol\theta_k$ is Pareto critical. This proves (ii).

We now prove (iii). By the descent lemma and assumption (A4), for every $k\ge 0$ we have
\begin{equation}
\Phi(\boldsymbol\theta_{k+1})
\le
\Phi(\boldsymbol\theta_k)
+
h\langle \nabla \Phi(\boldsymbol\theta_k),\mathbf d_k\rangle
+
\frac{Lh^2}{2}\|\mathbf d_k\|^2
\label{eq:63}
\end{equation}
From \eqref{eq:54} it follows that
\begin{equation}
\langle \nabla \Phi(\boldsymbol\theta_k),\mathbf d_k\rangle
\le
-\frac{1}{2}\|\mathbf d_k\|^2
\label{eq:64}
\end{equation}
Substituting \eqref{eq:64} into \eqref{eq:63} yields
\begin{equation}
\Phi(\boldsymbol\theta_{k+1})
\le
\Phi(\boldsymbol\theta_k)
-
\frac{h}{2}\|\mathbf d_k\|^2
+
\frac{Lh^2}{2}\|\mathbf d_k\|^2
\label{eq:65}
\end{equation}
Therefore,
\begin{equation}
\Phi(\boldsymbol\theta_{k+1})
\le
\Phi(\boldsymbol\theta_k)
-
\frac{h}{2}(1-Lh)\|\mathbf d_k\|^2
\label{eq:66}
\end{equation}
Since $h\le 1/(2L)$ by \eqref{eq:51}, we have
\begin{equation}
1-Lh\ge \frac{1}{2}
\label{eq:67}
\end{equation}
Combining \eqref{eq:66} and \eqref{eq:67} gives
\begin{equation}
\Phi(\boldsymbol\theta_{k+1})
\le
\Phi(\boldsymbol\theta_k)
-
\frac{h}{4}\|\mathbf d_k\|^2
\label{eq:68}
\end{equation}
Using \eqref{eq:55}, we further obtain
\begin{equation}
\Phi(\boldsymbol\theta_{k+1})
\le
\Phi(\boldsymbol\theta_k)
-
\frac{hM^2}{4}\|\nabla \Phi(\boldsymbol\theta_k)\|^2
\label{eq:69}
\end{equation}
Summing \eqref{eq:69} over $k=0,1,\ldots,K$ yields
\begin{equation}
\frac{hM^2}{4}\sum_{k=0}^K \|\nabla \Phi(\boldsymbol\theta_k)\|^2
\le
\Phi(\boldsymbol\theta_0)-\Phi(\boldsymbol\theta_{K+1})
\label{eq:70}
\end{equation}
By item (i), we have $\boldsymbol\theta_{K+1}\in\mathcal X$. Hence, by \eqref{eq:53},
\begin{equation}
\Phi(\boldsymbol\theta_{K+1})\ge \Phi_*
\label{eq:71}
\end{equation}
Substituting \eqref{eq:71} into \eqref{eq:70} gives
\begin{equation}
\frac{hM^2}{4}\sum_{k=0}^K \|\nabla \Phi(\boldsymbol\theta_k)\|^2
\le
\Phi(\boldsymbol\theta_0)-\Phi_*
\label{eq:72}
\end{equation}
Dividing both sides of \eqref{eq:72} by $K+1$ proves \eqref{eq:56}.

Finally, since
\begin{equation}
\min_{0\le k\le K}\|\nabla \Phi(\boldsymbol\theta_k)\|^2
\le
\frac{1}{K+1}\sum_{k=0}^K \|\nabla \Phi(\boldsymbol\theta_k)\|^2
\label{eq:73}
\end{equation}
combining \eqref{eq:73} with \eqref{eq:56} yields
\begin{equation}
\min_{0\le k\le K}\|\nabla \Phi(\boldsymbol\theta_k)\|^2
\le
\frac{4\bigl(\Phi(\boldsymbol\theta_0)-\Phi_*\bigr)}{hM^2(K+1)}
\label{eq:74}
\end{equation}
Taking square roots on both sides of \eqref{eq:74} proves \eqref{eq:57}.
\end{proof}

\begin{remark}
\textcolor{blue}{\textbf{Algorithm}}~\ref{algorithm:2} is the exact step-search rule used to define an idealized monotone update. It selects the largest step size such that, along the entire segment \([\,\boldsymbol\theta_k,\boldsymbol\theta_k+h\mathbf d_k\,]\), the directional derivatives of all objectives remain nonpositive. When all objective functions \(\mathcal F_t\) are convex and differentiable, this rule can be simplified. Indeed, for
\[
\varphi_{k,t}(\tau)\triangleq \mathcal F_t(\boldsymbol\theta_k+\tau \mathbf d_k),
\]
convexity of \(\mathcal F_t\) implies that \(\varphi_{k,t}\) is convex in \(\tau\), and therefore \(\varphi_{k,t}'(\tau)\) is nondecreasing on its domain. Hence, if for some \(h\in[0,1]\),
\[
\left\langle \nabla \mathcal F_t(\boldsymbol\theta_k+h\mathbf d_k),\mathbf d_k\right\rangle \le 0
\qquad t=1,\ldots,T,
\]
then for every \(\tau\in[0,h]\),
\[
\left\langle \nabla \mathcal F_t(\boldsymbol\theta_k+\tau\mathbf d_k),\mathbf d_k\right\rangle
=
\varphi_{k,t}'(\tau)
\le
\varphi_{k,t}'(h)
=
\left\langle \nabla \mathcal F_t(\boldsymbol\theta_k+h\mathbf d_k),\mathbf d_k\right\rangle
\le 0.
\]
Therefore, in the convex case, it is sufficient to check only the endpoint condition. This justifies a simpler backtracking implementation that starts from \(h=1\) and repeatedly shrinks \(h\) until the endpoint test is satisfied. By contrast, when the objectives are nonconvex, the endpoint test is no longer sufficient, and one must revert to the interval-based condition in \textcolor{blue}{\textbf{Algorithm}}~\ref{algorithm:2}, or to a numerical approximation of that condition.

This distinction is consistent with the implementations used in this paper. For deterministic convex applications, such as the multi-criteria traffic assignment problem in \textcolor{blue}{\textbf{Section}}~\ref{sec:5.2}, we use a practical backtracking step-size rule based only on the endpoint directional derivatives. For more general nonconvex cases, the exact interval condition may instead be approximated numerically. On the other hand, in stochastic gradient-based training problems, such as the physics-informed machine-learning-based car-following model in \textcolor{blue}{\textbf{Section}}~\ref{sec:5.1}, we do not perform an exact line search. After computing a common descent direction, the trainable parameters are updated by an adaptive optimizer, so the effective step size is determined by the optimizer rather than by \textcolor{blue}{\textbf{Algorithm}}~\ref{algorithm:2}. Accordingly, \textcolor{blue}{\textbf{Algorithm}}~\ref{algorithm:2} should be viewed as an idealized exact step-search rule for deterministic monotone updates, whereas \textcolor{blue}{\textbf{Theorem}}~\ref{theorem:2} provides a convergence-rate result for a fixed-step theoretical variant of the proposed update.
\end{remark}

\par  \textcolor{blue}{\textbf{Algorithm}} \ref{algorithm:1} summarizes the overall min-max progamming framework.  In each iteration, it first solves a min–max or quadratic programming problem to determine the direction $\mathbf{d}$ that best decreases all objectives while satisfying the hard constraints.  This iterative process continues until convergence criteria are met: either a Pareto-critical point is reached, or the algorithm exceeds the maximum number of iterations. It is important to note that our algorithms update the current solution, denoted as \(\boldsymbol{\theta}\), by moving along the direction \(\mathbf{d}\) with a step size \(h\). This approach ensures a monotonic decrease for all objectives over the interval \((0, h)\) (as indicated in \cite{desideri2012multiple} and shown in equation \((\ref{eq:11})\)), which can be implemented though an exact search algorithm shown in \textcolor{blue}{\textbf{Algorithm}} \ref{algorithm:2}\footnote{\textcolor{blue}{\textbf{Algorithm}}~\ref{algorithm:2} is the exact step-search rule used in the convergence analysis. 
In numerical implementation, one may introduce a small tolerance $\epsilon>0$ and terminate the line search when the computed step size is smaller than $\epsilon$. 
Such a truncation is purely for numerical robustness.}. It is well known that a class of $\min-\max$ programs can be equivalently reformulated as a single-level constrained optimization problem by introducing auxiliary variables and a set of linear constraints. Specifically, the inner maximization can be replaced by a set of linear inequalities that characterize its optimality conditions or feasible region, thereby converting the original bilevel structure into a tractable single-level program. In our case, by introducing the auxiliary variable $\eta$, \textcolor{blue}{\textbf{Problem}} \eqref{eq:6} can be reformulated as \textcolor{blue}{\textbf{Problem}} \eqref{eq:75}, which is a pure linear program. 

\begin{algorithm}
    \caption{\textsc{Min-Max programming}}
    \label{algorithm:1}
    \LinesNumbered
    \KwIn {$\boldsymbol{\theta}_{0}$,Gradient vectors $\mathbf{P} = [\mathbf{p}_{1}, \mathbf{p}_{2},...,\mathbf{p}_{T}]|_{\boldsymbol{\theta} = \boldsymbol{\theta}_{0}} = [\nabla \mathcal{F}^1(\boldsymbol{\theta}),\nabla \mathcal{F}^2(\boldsymbol{\theta}),...,\nabla \mathcal{F}^T(\boldsymbol{\theta})]|_{\boldsymbol{\theta} = \boldsymbol{\theta}_{0}}$, $K$,$tol$}{
    $k = 0, \boldsymbol{\theta} = \boldsymbol{\theta}_{0}$ \\
    \While{$k \leq M_1$}{
    Solve: \\
\begin{equation}
\begin{aligned}
w = \min_{\eta,\, \mathbf{d}} \quad & \eta \\
\text{s.t.} \quad & f_{i}(\boldsymbol{\theta}+\mathbf{d})  \leq 0, \quad i = 1,2,...,m\\
& \nabla \mathcal{F}^t(\boldsymbol{\theta})^\top \mathbf{d} \leq \eta, \quad \forall t, \\
& \nabla \mathcal{F}^t(\boldsymbol{\theta})^\top \mathbf{d} \leq 0, \quad \forall t, \\
& \mathbf{-1} \leq \mathbf{d}  \leq \mathbf{1}
\end{aligned}
\label{eq:75}
\end{equation} \\
\If{$|w| \leq tol$}{\textbf{break}}
$h = $\textsc{Step search algorithm}($\mathbf{d^{*}}$) \\
\If{$h = 0$}{
\textbf{break} \\
}
$\boldsymbol{\theta} = \boldsymbol{\theta} + h \cdot \mathbf{d}^{*} $ 
}
}
\KwOut{$\boldsymbol{\theta}$}
\end{algorithm}

\begin{algorithm}[t]
    \caption{\textsc{Step search algorithm}}
    \label{algorithm:2}
    \LinesNumbered
    \KwIn{$\mathbf d^{*}$,$\boldsymbol{\theta}$,$\mathcal{F}_{t},t = 1,...,T$}
    
    Compute
    \begin{equation}
    h
    =
    \max \Bigl\{
    h' \in [0,1]\ \Big|\ 
    \nabla \mathcal{F}_{t}(\boldsymbol{\theta} + \tau \mathbf d^{*})^{\top}\mathbf d^{*} \le 0,\ 
    \forall \tau \in [0,h'],\ t = 1,2,\ldots,T
    \Bigr\}
    \label{eq:76}
    \end{equation}
    
    \KwOut{$h$}
\end{algorithm}

\section{Computational Efficiency of the LP-Based Min-max Algorithm} \label{sec:4}
\par In this section, we analyze the computational complexity of solving the min-max programming in \textcolor{blue}{\textbf{Problem}}~\eqref{eq:6} and compare it with the classical steepest descent formulation in \textcolor{blue}{\textbf{Problem}}~\eqref{eq:1} proposed by \cite{fliege2000steepest}. Since large-scale linear programs and convex quadratic programs are commonly solved by primal-dual interior-point methods, both the theoretical complexity discussion and the numerical timing analysis in this section are based on the primal-dual interior-point method.

\par Multi-objective traffic assignment problems often impose linear constraints with highly sparse matrices. Let
\begin{equation}
\mathbf u
\triangleq
\begin{bmatrix}
\mathbf d \\[2pt]
\eta
\end{bmatrix}
\in \mathbb R^{\,n+1}
\label{eq:77}
\end{equation}
Then \textcolor{blue}{\textbf{Problem}}~\eqref{eq:6} can be written in the compact LP inequality form
\begin{equation}
\begin{aligned}
\min_{\mathbf u}\quad & c^\top \mathbf u \\
\textnormal{s.t.}\quad & B\mathbf u \le q
\end{aligned}
\label{eq:78}
\end{equation}
where
\begin{equation}
c
\triangleq
\begin{bmatrix}
\mathbf 0_n \\[2pt]
1
\end{bmatrix}
\in\mathbb R^{\,n+1}
\qquad
B\in\mathbb R^{\,M\times(n+1)}
\qquad
q\in\mathbb R^{\,M}
\qquad
M:=m+2T+2n
\label{eq:79}
\end{equation}
Here, the \(m\) rows of \(B\mathbf u\le q\) correspond to the feasibility constraints, the \(2T\) rows correspond to the directional derivative constraints, and the \(2n\) rows correspond to the box constraints on \(\mathbf d\). The vector \(\mathbf u\) is unrestricted in sign.

\par Introducing primal slack variables \(s\in\mathbb R^{\,M}\) with \(s\ge 0\), \eqref{eq:78} is equivalently written as
\begin{equation}
\begin{aligned}
\min_{\mathbf u,s}\quad & c^\top \mathbf u \\
\textnormal{s.t.}\quad & B\mathbf u + s = q \\
& s \ge 0
\end{aligned}
\label{eq:80}
\end{equation}

\par Since \(\mathbf u\) is unrestricted in sign, the dual problem of \eqref{eq:80} is
\begin{equation}
\begin{aligned}
\max_{\lambda}\quad & -q^\top \lambda \\
\textnormal{s.t.}\quad & B^\top \lambda + c = 0 \\
& \lambda \ge 0
\end{aligned}
\label{eq:81}
\end{equation}
where \(\lambda\in\mathbb R^{\,M}\) is the Lagrange multiplier associated with the equality constraint \(B\mathbf u+s=q\).

\par The Karush--Kuhn--Tucker conditions for \eqref{eq:80} and \eqref{eq:81} are
\begin{subequations}
\begin{align}
& B^\top \lambda + c = 0
\label{eq:82a}\\
& B\mathbf u + s - q = 0
\label{eq:82b}\\
& s_i\lambda_i = 0
\qquad i=1,\ldots,M
\label{eq:82c}\\
& s\ge 0
\qquad
\lambda\ge 0
\label{eq:82d}
\end{align}
\label{eq:82}
\end{subequations}

\par For the primal-dual interior-point method, the perturbed complementarity conditions are
\begin{equation}
S\Lambda e = \sigma \mu e
\label{eq:83}
\end{equation}
where
\begin{equation}
S=\mathrm{diag}(s_1,\ldots,s_M)
\qquad
\Lambda=\mathrm{diag}(\lambda_1,\ldots,\lambda_M)
\qquad
e=(1,\ldots,1)^\top\in\mathbb R^{\,M}
\label{eq:84}
\end{equation}
and \(\sigma\in(0,1)\) is the centering parameter and \(\mu\) is the barrier parameter.

\par Define the residuals
\begin{equation}
r_d = B^\top\lambda + c
\qquad
r_p = B\mathbf u + s - q
\qquad
r_c = S\Lambda e - \sigma\mu e
\label{eq:85}
\end{equation}
Applying Newton's method to the perturbed KKT system yields the linear system
\begin{equation}
\begin{bmatrix}
0 & 0 & B^\top \\
B & I_M & 0 \\
0 & \Lambda & S
\end{bmatrix}_{(n+1+2M)\times(n+1+2M)}
\begin{bmatrix}
\Delta \mathbf u \\
\Delta s \\
\Delta \lambda
\end{bmatrix}
=
-
\begin{bmatrix}
r_d \\
r_p \\
r_c
\end{bmatrix}
\label{eq:86}
\end{equation}

\par We next compare \eqref{eq:86} with the classical steepest descent formulation in \textcolor{blue}{\textbf{Problem}}~\eqref{eq:1}. Its QP subproblem can be written as
\begin{equation}
\begin{aligned}
\min_{x\in\mathbb R^{\,n+1}}\quad
& \frac{1}{2}x^\top Gx + \hat c^\top x \\
\textnormal{s.t.}\quad
& \hat A x \ge \hat b
\end{aligned}
\label{eq:87}
\end{equation}
where
\begin{equation}
x
\triangleq
\begin{bmatrix}
\mathbf d \\[2pt]
\gamma
\end{bmatrix}
\in\mathbb R^{\,n+1}
\qquad
G
\triangleq
\begin{bmatrix}
I_n & 0 \\
0 & 0
\end{bmatrix}
\in\mathbb R^{\,(n+1)\times(n+1)}
\qquad
\hat c
\triangleq
\begin{bmatrix}
\mathbf 0_n \\[2pt]
1
\end{bmatrix}
\in\mathbb R^{\,n+1}
\label{eq:88}
\end{equation}
and
\begin{equation}
\hat A\in\mathbb R^{\,\hat M\times(n+1)}
\qquad
\hat b\in\mathbb R^{\,\hat M}
\qquad
\hat M:=m+T
\label{eq:89}
\end{equation}

\par Introducing primal slacks \(y\in\mathbb R^{\,\hat M}\) with \(y\ge 0\), the inequality constraint in \eqref{eq:87} becomes
\begin{equation}
\hat A x - y = \hat b
\label{eq:90}
\end{equation}
Let
\begin{equation}
Y=\mathrm{diag}(y_1,\ldots,y_{\hat M})
\qquad
\Lambda_q=\mathrm{diag}((\lambda_q)_1,\ldots,(\lambda_q)_{\hat M})
\qquad
e=(1,\ldots,1)^\top\in\mathbb R^{\,\hat M}
\label{eq:91}
\end{equation}
Then the corresponding perturbed KKT system can be written as
\begin{equation}
\begin{bmatrix}
G & 0 & -\hat A^\top \\
\hat A & -I_{\hat M} & 0 \\
0 & \Lambda_q & Y
\end{bmatrix}_{(n+1+2\hat M)\times(n+1+2\hat M)}
\begin{bmatrix}
\Delta x \\
\Delta y \\
\Delta \lambda_q
\end{bmatrix}
=
-
\begin{bmatrix}
r_d^{\,q} \\
r_p^{\,q} \\
Y\Lambda_q e-\sigma\mu e
\end{bmatrix}
\label{eq:92}
\end{equation}
where
\begin{equation}
r_d^{\,q}=Gx-\hat A^\top\lambda_q+\hat c
\qquad
r_p^{\,q}=\hat A x-y-\hat b
\label{eq:93}
\end{equation}

\par Therefore, the LP formulation \eqref{eq:78} uses
\begin{equation}
M=m+2T+2n
\label{eq:94}
\end{equation}
constraints, whereas the classical QP formulation \eqref{eq:87} uses
\begin{equation}
\hat M=m+T
\label{eq:95}
\end{equation}
constraints. Hence the proposed LP formulation introduces
\begin{equation}
M-\hat M=T+2n
\label{eq:96}
\end{equation}
additional constraints relative to the classical QP formulation.

\par In a primal-dual interior-point method, the computation of the Newton step, that is, solving the associated KKT linear system such as \eqref{eq:86} or \eqref{eq:92}, typically dominates the per-iteration computational cost. For the full KKT systems arising here, a more standard numerical treatment is to use symmetric indefinite factorization or, after suitable elimination, a Schur-complement-based linear-algebra strategy, rather than modified Cholesky factorization. In this sense, the classical QP-based steepest-descent subproblem has a simpler Newton system than the proposed LP formulation, because it involves fewer explicit constraints and therefore a smaller algebraic system.

\par The KKT system \eqref{eq:86} has dimension
\begin{equation}
n+1+2M
=
n+1+2(m+2T+2n)
\label{eq:97}
\end{equation}
whereas the KKT system \eqref{eq:92} has dimension
\begin{equation}
n+1+2\hat M
=
n+1+2(m+T)
\label{eq:98}
\end{equation}
Since \(T\ll n\) and \(m\asymp n\) in typical multi-criteria traffic assignment problems, the LP-based formulation yields a substantially larger KKT system than the classical QP-based formulation.

\par From a theoretical worst-case complexity perspective, for dense instances, the arithmetic complexity of an interior-point method applied to a linear program with \(p\) variables and \(k\) constraints is bounded by~\cite{nemirovski2004interior}
\begin{equation}
\mathcal O\!\left(
\sqrt{k}\,
\log\!\left(\frac{1}{\varepsilon}\right)
\cdot
\left(
p^2k + pk^2 + k^3
\right)
\right)
\label{eq:99}
\end{equation}
where \(\varepsilon\) is the target solution accuracy.

\par In our LP formulation,
\begin{equation}
p=n+1
\qquad
k=M=m+2T+2n
\label{eq:100}
\end{equation}
Hence the worst-case arithmetic complexity of solving \textcolor{blue}{\textbf{Problem}}~\eqref{eq:6} is
\begin{equation}
\mathcal O\!\left(
\sqrt{M}\,
\log\!\left(\frac{1}{\varepsilon}\right)
\cdot
\left(
(n+1)^2M + (n+1)M^2 + M^3
\right)
\right)
\label{eq:101}
\end{equation}
Assuming \(T\ll n\) and \(m\asymp n\), we have \(M=\Theta(n)\), so \eqref{eq:101} simplifies to
\begin{equation}
\mathcal O\!\left(
n^{3.5}\log\!\left(\frac{1}{\varepsilon}\right)
\right)
\label{eq:102}
\end{equation}

\par Under the same dense model, the classical QP formulation \eqref{eq:87} exhibits the same polynomial order in \(n\), since \(\hat M=m+T=\Theta(n)\) and its Newton system has comparable cubic linear-algebra scaling. Thus, the LP and QP formulations have the same dense worst-case asymptotic complexity order.

\par It should also be noted that many textbook presentations assume full-row-rank constraint matrices. In our setting, however, the number of constraints may substantially exceed the number of variables, so neither \(B\) nor \(\hat A\) needs to have full row rank. Modern solvers therefore typically apply presolve or preprocessing steps, such as removing redundant or linearly dependent constraints, before invoking the interior-point iterations. Since the exact reduction depends strongly on the particular problem instance, the comparison above is based on the original formulations \eqref{eq:78} and \eqref{eq:87}. In particular, while both formulations share the \(m\) feasibility constraints, preprocessing need not reduce the two models in the same proportion. Nevertheless, the LP formulation contains \(T+2n\) additional explicit constraints that are absent from the classical QP formulation, so its associated Newton system is generally larger in algebraic dimension.

\par Therefore, the proposed LP reformulation does not yield an inherent worst-case complexity advantage over the classical QP formulation under the dense model considered here. Rather, its advantage lies in the robustness and balance of the common descent direction that it generates. At the same time, when the problem size is small to moderate and no additional exploitable structure is used, the classical QP formulation often enjoys lower wall-clock time because of its smaller KKT systems and lower linear-algebra setup overhead. This statement should be interpreted as a structural tendency under the dense worst-case model, not as a universal ranking for all instances.

\par On the other hand, for the large-scale transportation problems that motivate this study, the constraint matrices are typically highly sparse and structured. In such settings, the practical performance of interior-point methods depends strongly on the sparsity pattern, fill-in behavior, and eliminable block structure of the associated Newton systems. Consequently, the dense worst-case comparison above should be interpreted only as a baseline theoretical benchmark. The empirical performance comparison must instead be assessed through the numerical experiments reported in \textcolor{blue}{\textbf{Section}}~\ref{sec:5}.

\section{Numerical examples} \label{sec:5}
\par In this section, we present two high-dimensional transportation numerical examples formulated as unconstrained and constrained multi-objective optimization problems: training of a physics-informed machine learning-based car-following model, and a multi-criteria traffic assignment problem. 
\subsection{Training of a physics-informed machine learning-based car-following model} \label{sec:5.1}

\begin{figure}
    \centering
    \includegraphics[width=0.9\linewidth]{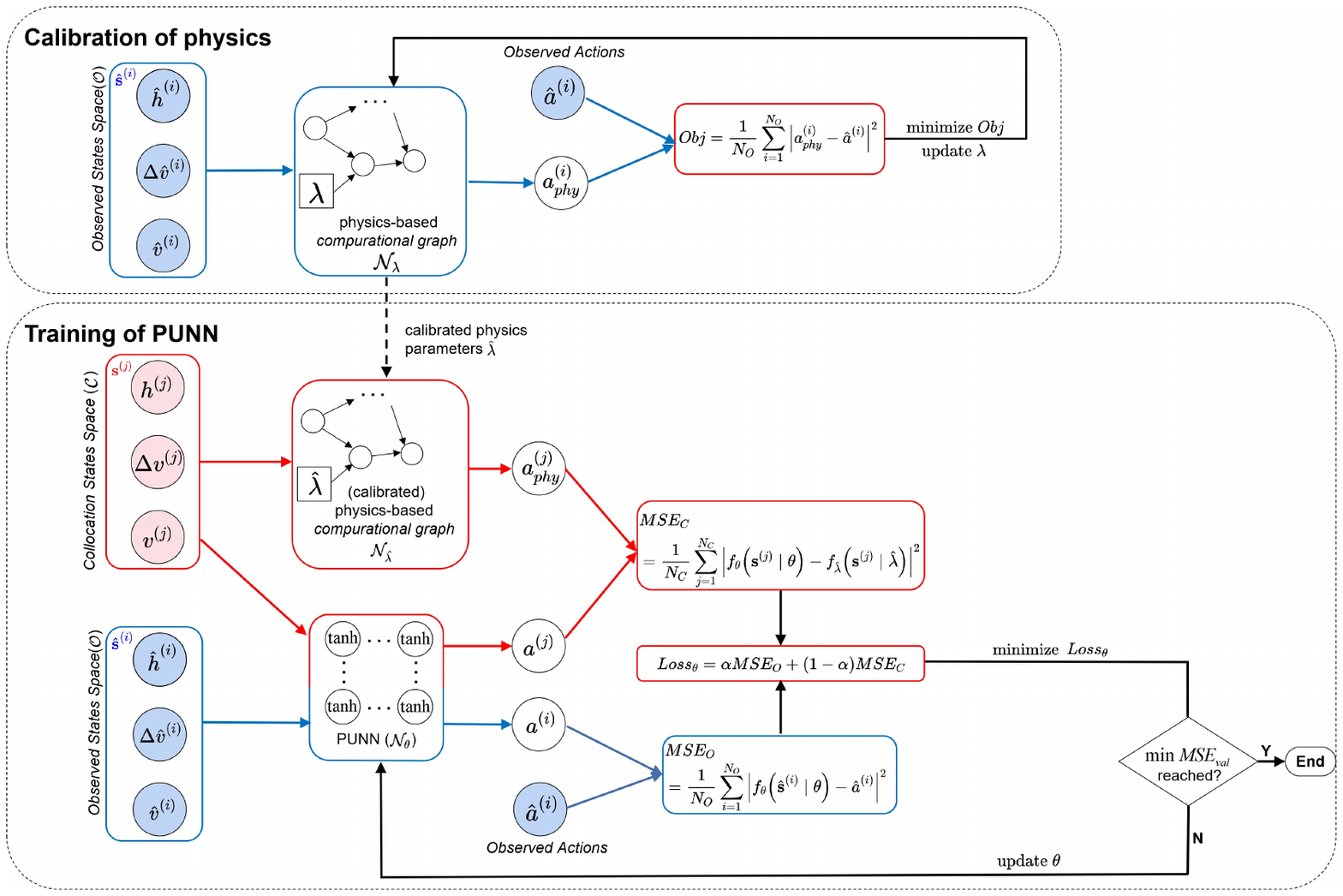}
    \caption{PINN-PICF architecture, copy from (\cite{mo2021physics})}
    \label{fig:1}
\end{figure}

As the first application of the proposed algorithm, we consider the training of a physics-informed machine learning (PIML) car-following model based on the framework of \cite{mo2021physics}. Following the viewpoint of \cite{sener2018multi}, the training of a multi-task model can be interpreted as a multi-objective optimization problem, where each task-specific loss is treated as an objective and the parameter is updated along a common descent direction in each iteration. In the present PIML setting, the two objectives correspond to the data-fitting term and the physics-based regularization term. For a generic PIML car-following model, the training loss is usually written as
\begin{equation}
\mathcal{L}(\boldsymbol{\theta})
=
\alpha\,\mathcal{L}_{\mathrm{data}}(\mathbf{X};\boldsymbol{\theta})
+
\beta\,\mathcal{L}_{\mathrm{physics}}(\mathbf{X};\boldsymbol{\theta})
\label{eq:103}
\end{equation}
where \(\boldsymbol{\theta}\) denotes the trainable parameters, \(\mathcal{L}_{\mathrm{data}}(\mathbf{X};\boldsymbol{\theta})\) is the data-fitting term, \(\mathcal{L}_{\mathrm{physics}}(\mathbf{X};\boldsymbol{\theta})\) is the physics-based term, and \(\alpha\) and \(\beta\) are scalar coefficients controlling their relative influence. Once MGDA-type methods are introduced, however, the update direction can be computed directly from the gradients of these two objectives, so that the training process no longer needs to rely on a pre-specified linear scalarization through fixed coefficients \(\alpha\) and \(\beta\). In this subsection, we adopt the open-source prediction-only mode in \cite{mo2021physics} as the baseline model\footnote{For baseline, the training process is essentially to optimize an unconstrained single objective optimization problem, in which the objective function is \eqref{eq:103}. The baseline can be found at \url{https://github.com/CU-DitecT/TRC21-PINN-CFM}.}, whose architecture is shown in \textcolor{blue}{\textbf{Figure}}~\ref{fig:1}, and retain the parameter settings provided in the released code as the reference configuration. On this basis, we incorporate several MGDA-based training strategies\footnote{For those training strategies based on different multi-gradient descent algorithms, the training process is essentially to optimize an unconstrained multi-objective optimization problem, where the two objectives are $\mathcal{L}_{\mathrm{data}}(\mathbf{X};\boldsymbol{\theta})$ and $\mathcal{L}_{\mathrm{physics}}(\mathbf{X};\boldsymbol{\theta})$. We use different multi-gradient descent algorithms to obtain a common descent direction $\mathbf{d}^{*}$ that can minimize all objectives in each iteration, and then update the model parameter through $\boldsymbol{\theta} = \boldsymbol{\theta} + h\cdot\mathbf{d}^{*}$ where $h$ is the step size.}, including the proposed min-max programming approach and classical MGDA variants, for comparative evaluation.
It should be emphasized that, unlike standard multi-task learning, the individual values of the data loss and physics loss are not the final quantities of interest in this application. What ultimately matters is the predictive performance of the trained car-following model on unseen trajectories. Therefore, we evaluate each training method on the test set using the root mean square errors (RMSE) of the predicted follower position and speed. Specifically, we utilize the Root Mean Square Error (RMSE) to measure the errors in both position and velocity as follows:
\begin{equation} \label{eq:104}
     RMSE_x = \sqrt{
            \frac{1}{N_T T_i}\sum_{i=1}^{N_T} \sum_{t=0}^{T_i}
            \left| 
                x(t)-\hat{x}(t)
            \right|^2 
     }
 \end{equation}   
 \begin{equation} \label{eq:105}
     RMSE_v = \sqrt{
            \frac{1}{N_T T_i}\sum_{i=1}^{N_T} \sum_{t=0}^{T_i}
            \left| 
                v(t)-\hat{v}(t)
            \right|^2 
     }
 \end{equation}
where \(x\) and \(v\) denote the positions and velocities of the follower, \(T_i\) is the time horizon of the \(i\)th trajectory, and \(N_T\) is the number of trajectories. We use the arithmetic mean of \(RMSE_x\) and \(RMSE_v\) as the final evaluation metric. For each training method, we conduct 30 independent runs, each corresponding to a different random seed. Moreover, although all experiments in this subsection use exactly the same dataset as \cite{mo2021physics}, our test set is constructed differently: instead of following the original split, we randomly extract 50 complete trajectories under a fixed random seed and use them consistently for all compared methods. The corresponding experimental results are summarized in \textcolor{blue}{\textbf{Table}}~\ref{table:2}, and all experiments conducted in this subsection is based on NVIDA GeForce RTX 4080 GPU.\footnote{Since our interest lies in the best achievable predictive performance of each training method under the same dataset and random-seed setting, rather than in its performance under a default hyperparameter choice, we tuned the main training hyperparameters for each method, including the learning rate, the maximum number of training iterations, and the patience parameter in early stopping, and report the best final result obtained under that fixed experimental setup over 30 runs.}
\begin{table}[htbp]
\centering
\caption{Comparison between different training methods} 
\label{table:2}
\begin{tabular}{ccc}
\thicktoprule
Training Method & Optimal coefficient & Average RSME \\ \hline
Baseline (\cite{mo2021physics}) & $\alpha = 0.9,\beta = 0.1$  &  4.9755 $\pm$ 1.6838\\
Ours & -  &  \textbf{4.6863 $\pm$ 0.8208}\\
Multiple-gradient descent algorithm (\cite{desideri2012multiple}) & -  &  6.4096 $\pm$ 6.2854 \\
The steepest descent algorithm (\cite{fliege2000steepest}) & -  & 5.2355 $\pm$ 1.9150 \\
% Dual cone gradient descent (Center) \cite{hwang2024dual} & -  &   \\
% Dual cone gradient descent (Average)\cite{hwang2024dual} & -  &   \\
% Dual cone gradient descent (Projection)\cite{hwang2024dual} & -  &   \\
\thickbottomrule
\end{tabular}
\end{table}
\begin{figure}[htbp]
    \centering
    \includegraphics[width=0.7\linewidth]{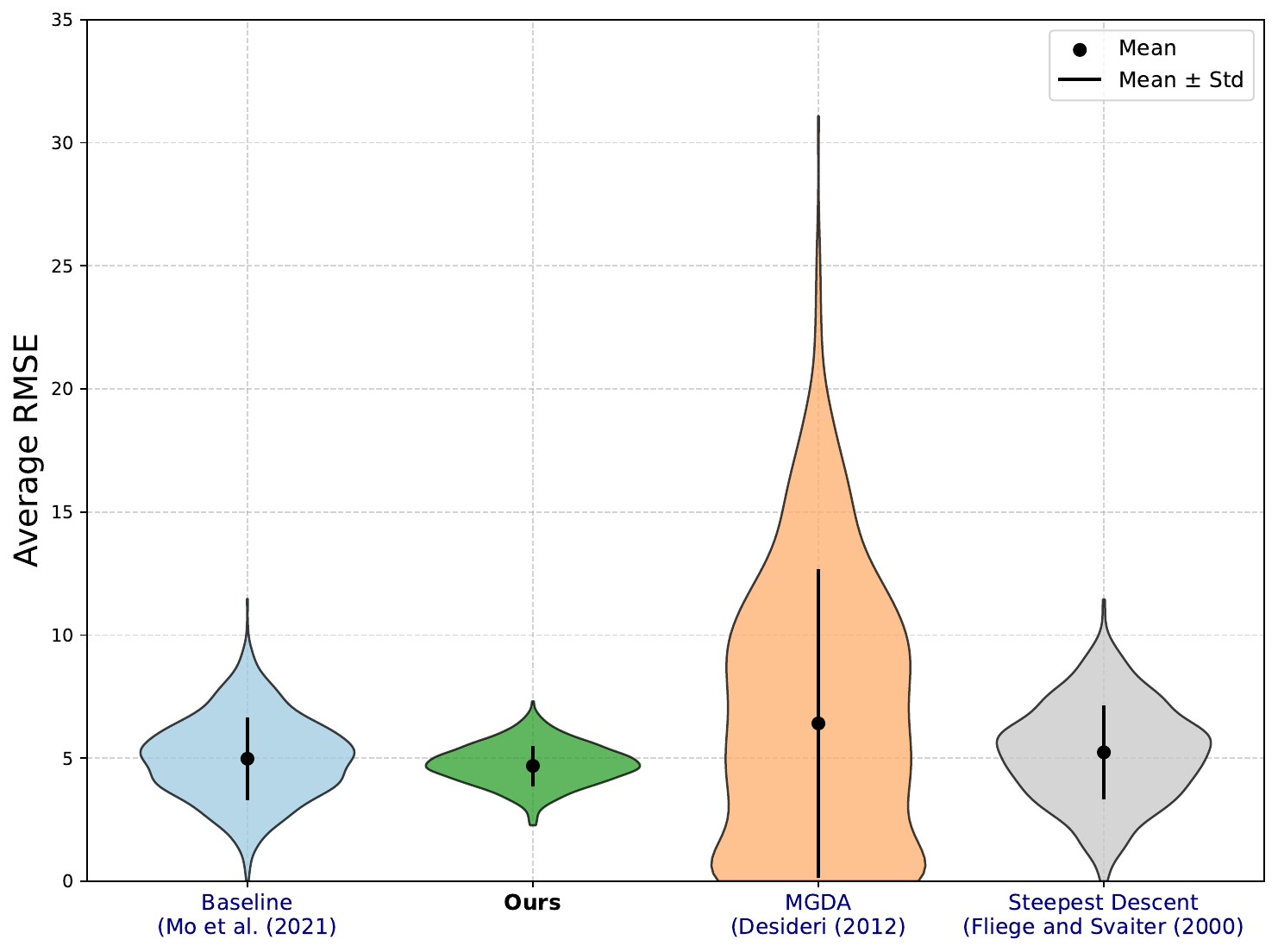}
    \caption{Comparison of the average RMSE distributions across different training methods. The proposed method achieves the lowest mean RMSE and exhibits smaller variability across random seeds compared with the baseline and the classical MGDA-based methods.}
    \label{fig:2}
\end{figure}
\par The results in \textcolor{blue}{\textbf{Table}}~\ref{table:2} and \textcolor{blue}{\textbf{Figure}}~\ref{fig:2} show that the proposed min--max programming method achieves the best average predictive performance among all compared training strategies. Compared with the baseline from \cite{mo2021physics}, whose best scalarization coefficients are \(\alpha=0.9\) and \(\beta=0.1\), our method reduces the average RMSE from \(4.9755\) to \(4.6863\), corresponding to a relative improvement of approximately \(5.81\%\). In addition, the standard deviation decreases from \(1.6838\) to \(0.8208\), indicating that the proposed method also produces more stable training outcomes across different random seeds. Compared with the classical multiple-gradient descent algorithm of \cite{desideri2012multiple} and the steepest descent algorithm of \cite{fliege2000steepest}, the relative reductions in average RMSE are approximately \(26.89\%\) and \(10.49\%\), respectively. \textcolor{blue}{\textbf{Figure}}~\ref{fig:2} further illustrates the empirical distribution of the repeated runs, highlighting that the proposed method achieves both the lowest average RMSE and a relatively concentrated performance distribution.

The two classical MGDA-based methods do not outperform the baseline in this experiment. In particular, the multiple-gradient descent algorithm of \cite{desideri2012multiple} produces a relatively large average RMSE and a high standard deviation, suggesting that its performance is sensitive to random initialization and training randomness. This should not be interpreted as evidence that classical MGDA methods are generally ineffective. Rather, it reflects a fundamental difference between deterministic multi-objective optimization and machine-learning-based training. Although the PIML training process can be formulated as a multi-objective optimization problem by treating the data-fitting loss and the physics-based loss as two objectives, these loss terms are constructed only from the training data. The final criterion of interest, however, is the trajectory-level prediction performance on the test set. Therefore, a direction that decreases the training losses does not necessarily guarantee better generalization performance or lower test RMSE.

The proposed min-max programming method appears to be more suitable for this setting because of the way it selects a common descent direction. Instead of only requiring a direction that decreases all objectives, it seeks a direction that makes the least improved objective decrease as much as possible by minimizing the maximum directional derivative. In other words, among all common descent directions, the proposed method favors a direction that improves the worst-case objective most aggressively. In this PIML car-following experiment, such a balanced worst-case descent mechanism leads to better test-set prediction performance and lower variability across random seeds. This result suggests that, when applying multi-objective optimization ideas to machine-learning training, the design of the common descent direction can have a substantial impact on the final predictive performance, even when all methods are based on the same underlying loss components.

\subsection{Multi-criteria traffic assignment problem} \label{sec:5.2}

\par The traffic assignment problem (TAP) is a classical and extensively studied topic in transportation science. It seeks to determine how origin-destination travel demand is allocated across a transportation network so that the resulting flow satisfies a certain behavioral equilibrium condition, such as user equilibrium or system optimum, under congestion-dependent travel costs. As a fundamental element of traffic modeling, TAP plays a critical role in analyzing and forecasting network performance under various demand scenarios.

\par Over the decades, a wide range of TAP variants have been developed to reflect different modeling assumptions and behavioral perspectives. These include the static user equilibrium (SUE) (\cite{nie2010class, nie2012note, wang2025entropy}, etc.), the system optimum (SO) (\cite{yang2004multi, yang2008existence, guo2009user, ma2014continuous}, etc.), the stochastic user equilibrium (STUE) (\cite{ye2022stochastic, xiao2019day, yang2001simultaneous, gentile2018new}, etc.), and dynamic traffic assignment (DTA) (\cite{zhu2000existence, han2015formulation, friesz2024dynamic, zhong2011dynamic}, etc.).

\par While these single-objective models have proven valuable in theory and practice, real-world transportation planning often involves multiple and potentially conflicting performance criteria. For instance, the UE condition reflects the behavior of individual travelers minimizing their own travel costs, whereas the SO minimizes the total travel cost across the network. In many cases, UE leads to inefficiencies from a system-wide perspective, while SO typically requires centralized coordination that may be infeasible or unpopular.

\par To address competing considerations in traffic assignment, researchers have proposed multi-criteria formulations that explicitly balance various objectives, such as UE and SO. These approaches provide planners with a range of trade-offs between individual travel performance and the collective efficiency of the transportation system. While many existing studies handle such problems by applying linear scalarization to convert the multi-objective formulation into a single-objective one, this approach presumes fixed, linear trade-offs among objectives and often lacks theoretical or empirical justification for the chosen scalarization weights.

\par In contrast, we retain the intrinsic multi-objective nature of the problem and solve it directly using the min--max programming framework proposed in this paper. By avoiding fixed scalarization weights, the proposed method computes a common descent direction from the UE and SO objectives themselves. This provides a direct way to balance individual travel behavior and system-level efficiency without imposing a pre-specified linear trade-off.

\par In the second numerical example, we apply the proposed min--max programming method to a bi-criteria TAP that balances UE and SO objectives. The well-studied Sioux Falls network is used. This network comprises 24 nodes, 76 links, 528 origin-destination (O-D) pairs, and 360600 demand units.

\par Let \( \mathcal{N} \) be the node set and \( \mathcal{A} \) be the link set, whose elements are ordered tuples of the form \( (n_{1}, n_{2}) \), with \( n_{1}, n_{2} \in \mathcal{N} \). Origins and destinations are denoted by \( \mathcal{R} \subset \mathcal{N} \) and \( \mathcal{S} \subset \mathcal{N} \), respectively. An O-D pair is represented as \( (r, s) \), where \( r \in \mathcal{R} \) and \( s \in \mathcal{S} \). The path set is defined as \( \mathcal{K} \).

\par We adopt the Bureau of Public Roads (BPR) function (\cite{us1964traffic}) to model the relationship between link flow and link travel time:
\begin{equation}
    t_a(x_a)
    =
    t_{a0}
    \left(
    1+\alpha\left(\frac{x_a}{c_a}\right)^{\beta}
    \right),
    \qquad \forall a \in \mathcal{A}
    \label{eq:106}
\end{equation}
where \(t_{a0}\) and \(c_a\) denote the free-flow travel time and capacity of link \(a\), respectively, and \(x_a\) is the flow on link \(a\). We set \(\alpha=0.15\) and \(\beta=4\). Then, the UE and SO objectives can be written as follows:
\begin{subequations}\label{eq:107}
\begin{align}
  & \min\ \mathcal{F}_{1} = \sum_{a}\int_{0}^{x_{a}}t_{a}(x)\,dx
    \label{eq:107a}\\
  &\text{s.t.}\quad
    \sum_{k}f_{k}^{rs} = q_{rs},
    \quad \forall r \in \mathcal{R},\ s \in \mathcal{S}
    \label{eq:107b}\\
  &\quad
    x_{a} = \sum_{r}\sum_{s}\sum_{k}\delta_{a,k}^{rs}f_{k}^{rs},
    \quad \forall a \in \mathcal{A}
    \label{eq:107c}\\
  &\quad
    f_{k}^{rs} \geq 0,
    \quad \forall k \in \mathcal{K},\ r \in \mathcal{R},\ s \in \mathcal{S}
    \label{eq:107d}\\
  &\quad
    x_{a} \geq 0,
    \quad \forall a\in \mathcal{A}
    \label{eq:107e}
\end{align}
\end{subequations}
\begin{subequations}\label{eq:108}
\begin{align}
  & \min\ \mathcal{F}_{2} = \sum_{a}x_{a}t_{a}(x_{a})
    \label{eq:108a}\\
  &\text{s.t.}\quad
    \sum_{k}f_{k}^{rs} = q_{rs},
    \quad \forall r \in \mathcal{R},\ s \in \mathcal{S}
    \label{eq:108b}\\
  &\quad
    x_{a} = \sum_{r}\sum_{s}\sum_{k}\delta_{a,k}^{rs}f_{k}^{rs},
    \quad \forall a \in \mathcal{A}
    \label{eq:108c}\\
  &\quad
    f_{k}^{rs} \geq 0,
    \quad \forall k \in \mathcal{K},\ r \in \mathcal{R},\ s \in \mathcal{S}
    \label{eq:108d}\\
  &\quad
    x_{a} \geq 0,
    \quad \forall a\in \mathcal{A}
    \label{eq:108e}
\end{align}
\end{subequations}
where \(f_{k}^{rs}\) is the flow on path \(k\) connecting the O-D pair \((r,s)\), \(q_{rs}\) is the demand between the O-D pair \((r,s)\), and \(\delta_{a,k}^{rs}\) is the link-path incidence indicator: \(\delta_{a,k}^{rs}=1\) if link \(a\) belongs to path \(k\) for O-D pair \((r,s)\), and \(\delta_{a,k}^{rs}=0\) otherwise.

\par
\begin{figure}
    \centering
    \includegraphics[width=1.0\linewidth]{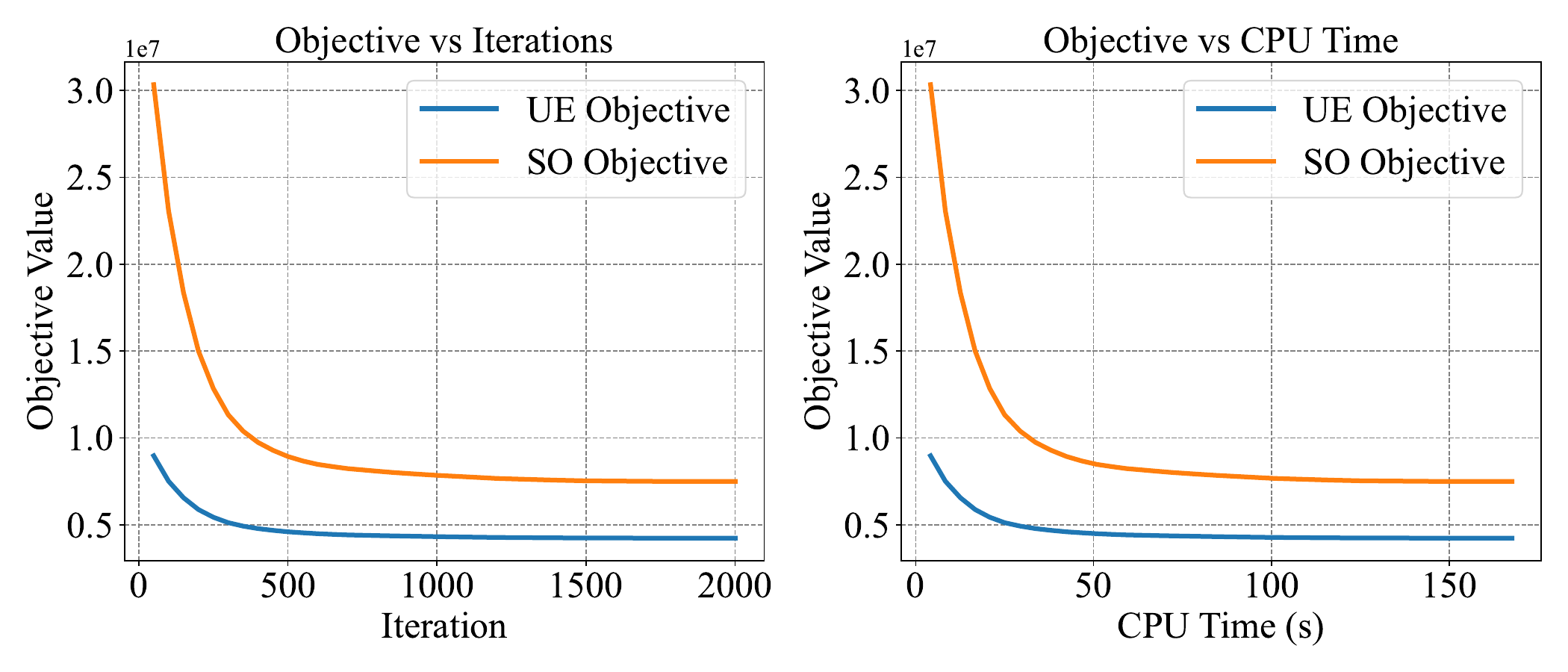}
    \caption{Comparison of the User Equilibrium (UE) and System Optimum (SO) objective values in the min--max multi-objective optimization: (i) plotted against iteration count; (ii) plotted against cumulative CPU time.}
    \label{figure:3}
\end{figure}

 To construct a set of candidate paths for each O-D pair, we implement a path enumeration algorithm based on the topology of the directed network. Starting from the link set and the O-D demand matrix, we build a directed graph with unit link weights. For each O-D pair with positive demand, we compute the shortest path using Dijkstra's algorithm and set a path-length threshold defined as a multiple of that shortest-path length, such as \(2.0\) times the shortest length. Subsequently, we employ Yen's algorithm (\cite{yen1970algorithm}) to generate up to \(K_{\mathrm{path}}=15\) shortest simple paths, ensuring that all enumerated paths are loopless and span from the origin to the destination. A total of 7920 paths are generated in \(\mathcal K\). As a result, at each iteration, the common descent direction for the multi-criteria traffic assignment problem is obtained by solving a linear program with 7921 decision variables and 24292 constraints. We initialize by assigning all demand to the shortest path. For this Sioux Falls example, we set the maximum number of iterations to \(N_{\mathrm{iter}}=2000\). \textcolor{blue}{\textbf{Figure}}~\ref{figure:3} compares the UE and SO objective values during the optimization process, plotted against iteration count and cumulative CPU time.

\par
\begin{figure}
    \centering
    \includegraphics[width=1.0\linewidth]{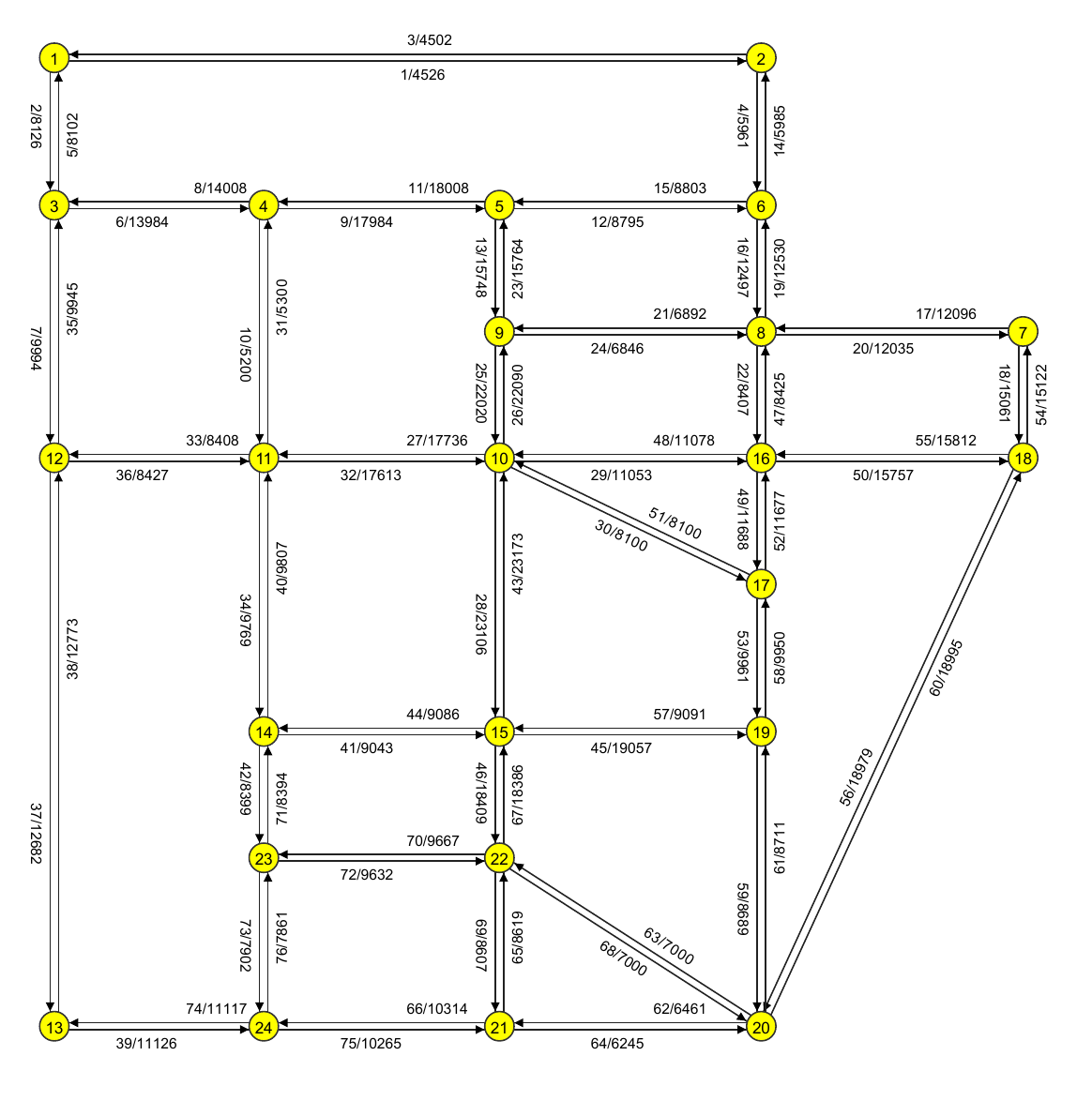}
    \caption{The final link flows on the Sioux Falls network.}
    \label{figure:4}
\end{figure}

 The experiment on the Sioux Falls network was run on an Apple M3 chip. The resulting link-flow pattern is shown in \textcolor{blue}{\textbf{Figure}}~\ref{figure:4}.

\subsection{Computational efficiency and robustness analysis}

\par To empirically validate the computational performance and robustness of the proposed LP-based formulation, we use the TAP as the representative application. We conducted numerical experiments on two relatively large-scale standard TAP instances: the Anaheim network (AN)\footnote{Refer to: \url{https://github.com/bstabler/TransportationNetworks/tree/master/Anaheim}} and the Eastern Massachusetts network (EMN)\footnote{Refer to: \url{https://github.com/bstabler/TransportationNetworks/tree/master/Eastern-Massachusetts}}. The Anaheim network consists of 416 nodes, 914 links, and a total trip demand of approximately 104694. The EMN network consists of 74 nodes, 258 links, and a total trip demand of approximately 65576. Similarly, for each O-D pair with positive demand, we generate up to 15 candidate paths. This yields direction-finding linear programs with 21091 decision variables and 64680 constraints for the AN network, and 16380 decision variables and 50254 constraints for the EMN network.

\par
\begin{figure}[htbp]
    \centering
    \includegraphics[width=0.88\linewidth]{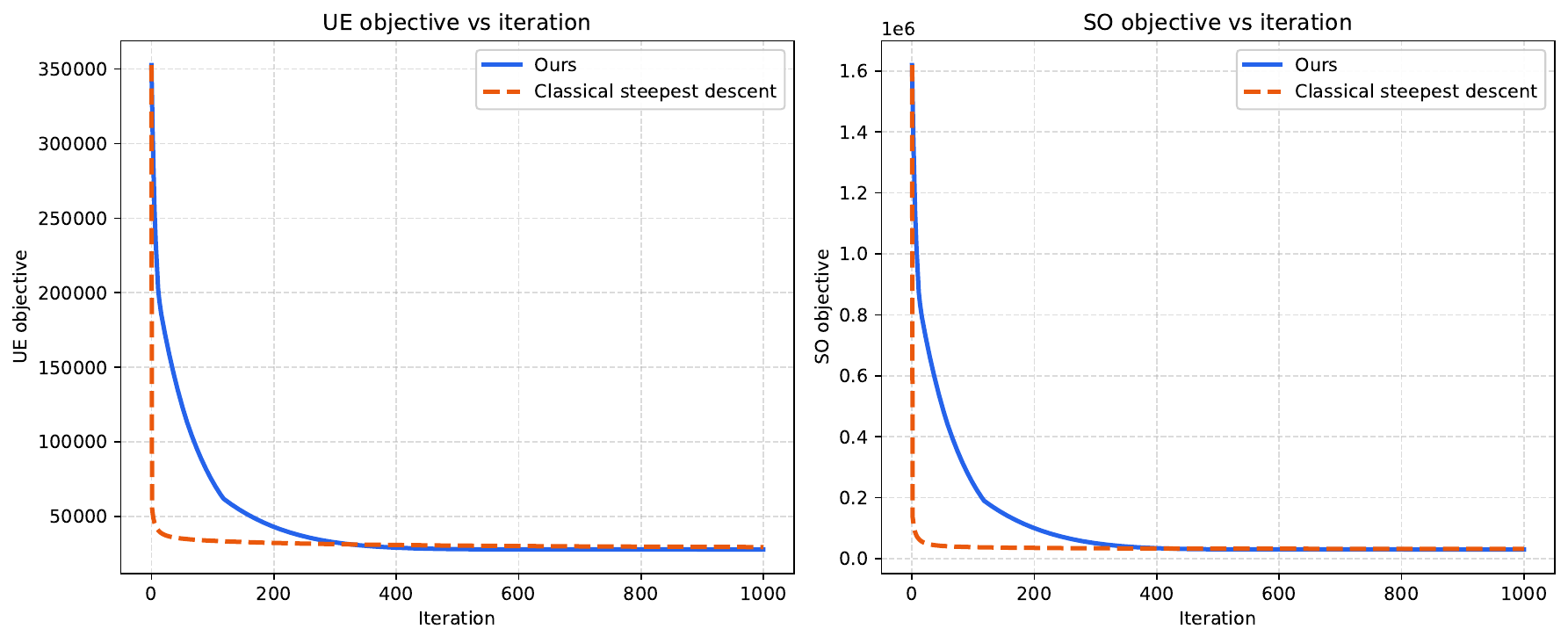}
    \caption{Convergence comparison between the proposed min--max programming method and the classical steepest descent algorithm on the Eastern Massachusetts network (Group I).}
    \label{fig:5}
\end{figure}

\begin{figure}[htbp]
    \centering
    \includegraphics[width=0.88\linewidth]{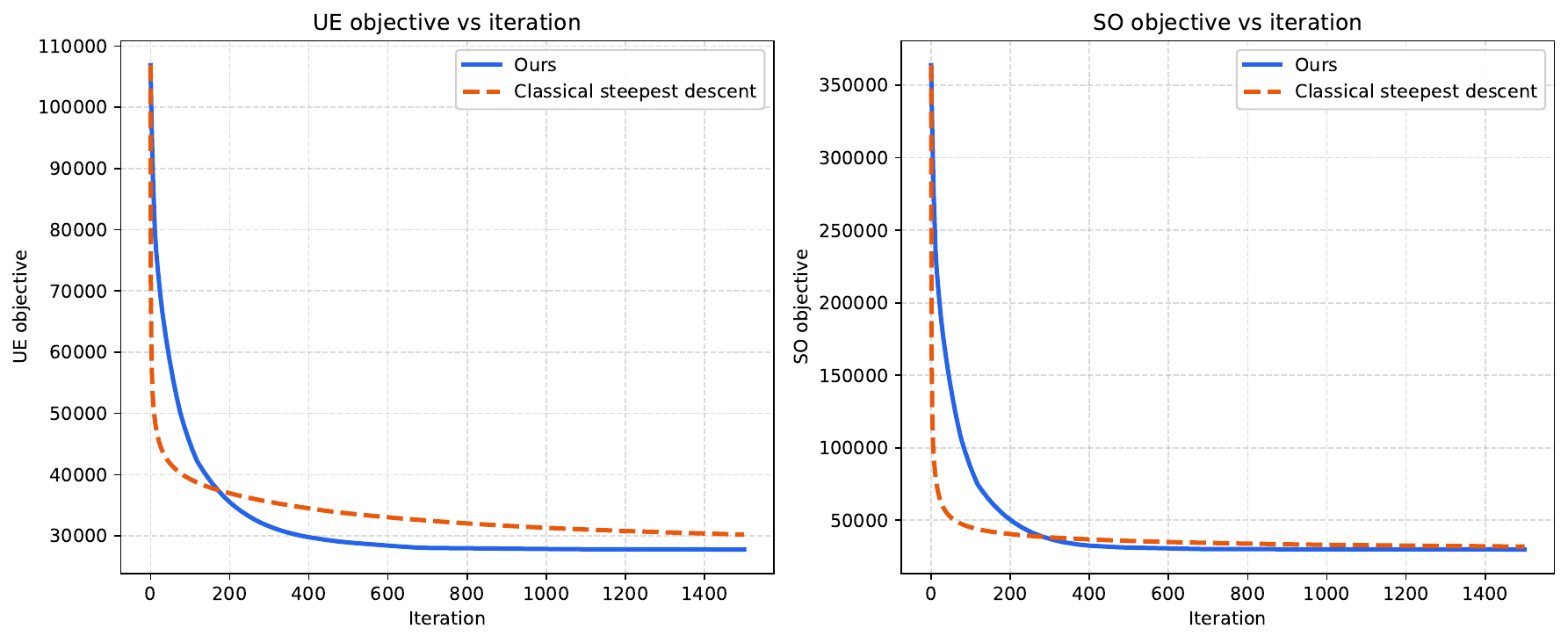}
    \caption{Convergence comparison between the proposed min--max programming method and the classical steepest descent algorithm on the Eastern Massachusetts network (Group II).}
    \label{fig:6}
\end{figure}

\begin{figure}[htbp]
    \centering
    \includegraphics[width=0.88\linewidth]{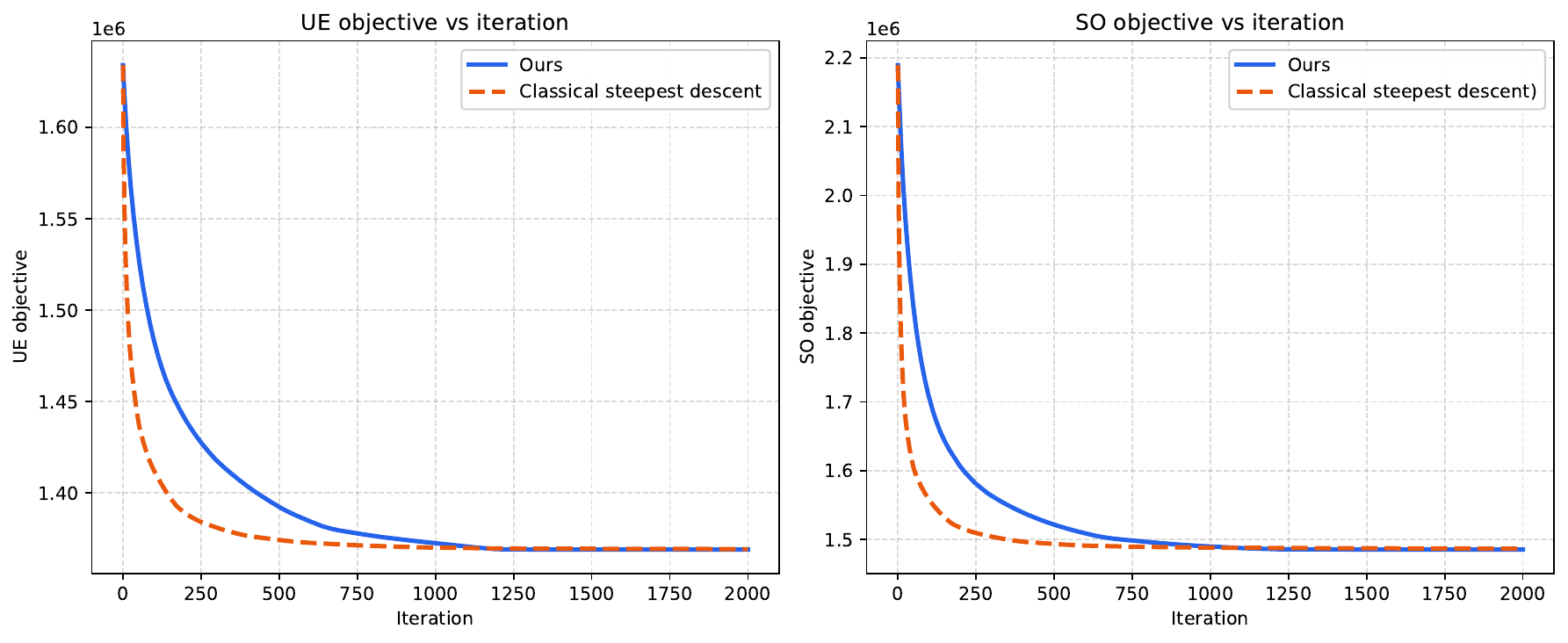}
    \caption{Convergence comparison between the proposed min--max programming method and the classical steepest descent algorithm on the Anaheim network (Group I).}
    \label{fig:7}
\end{figure}

\begin{figure}[htbp]
    \centering
    \includegraphics[width=0.88\linewidth]{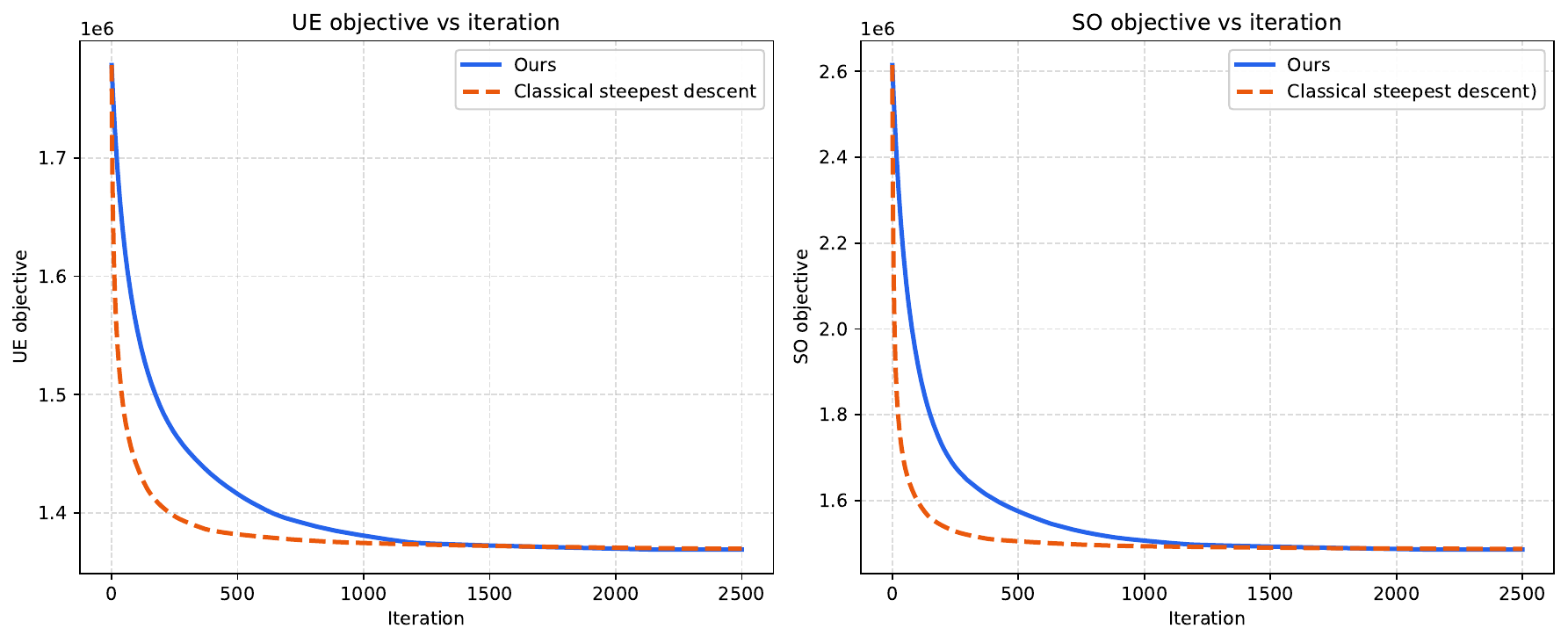}
    \caption{Convergence comparison between the proposed min--max programming method and the classical steepest descent algorithm on the Anaheim network (Group II).}
    \label{fig:8}
\end{figure}

 For each network, we consider two experimental groups. In Group I, we assign all O-D demand to the shortest-distance path under free-flow travel times. In Group II, we assign all O-D demand to the second shortest-distance path under free-flow travel times. If an O-D pair has no second-shortest-distance path, its demand is assigned to the shortest-distance path. As a baseline for comparison, we implement the classical steepest descent algorithm of \cite{fliege2000steepest}. All experiments were conducted on an Apple M3 chip, and the convergence tolerance is fixed across methods to ensure comparability. The numerical results are summarized in \textcolor{blue}{\textbf{Table}}~\ref{table:3} and \textcolor{blue}{\textbf{Table}}~\ref{table:4}. To further illustrate the dynamic convergence behavior of the compared methods, part of the convergence trajectories are shown in \textcolor{blue}{\textbf{Figures}}~\ref{fig:5}--\ref{fig:8}.

\par \textcolor{blue}{\textbf{Figures}}~\ref{fig:5}--\ref{fig:8} provide a detailed view of some of the convergence process behind the final objective values reported in \textcolor{blue}{\textbf{Table}}~\ref{table:3} and \textcolor{blue}{\textbf{Table}}~\ref{table:4}. A consistent trend can be observed across the four cases. In the early stage of the iteration process, the classical steepest descent algorithm often decreases the objective values more rapidly than the proposed min--max programming method. However, as the iteration proceeds, the proposed method generally catches up and then achieves better or more balanced final objective values. This empirical behavior is consistent with the design principle of the proposed direction-finding subproblem. Since the proposed method minimizes the maximum directional derivative among all objectives, it explicitly focuses on improving the least-improved objective at each iteration. As a result, the generated direction may be less aggressive in the initial stage, but it tends to produce more robust and balanced long-run progress across the competing UE and SO objectives.

\par This dynamic convergence pattern also helps explain the computational results discussed below. In particular, the faster initial decrease of the classical steepest descent algorithm is compatible with its lower per-iteration computational burden, whereas the proposed min--max programming method may require more effort per iteration but can generate more effective common descent directions and thus achieve better final convergence behavior under practical iteration budgets.

\par First, the results show that, for the same number of iterations, the classical steepest descent algorithm is generally faster in raw CPU time. This observation is consistent with the complexity discussion in \textcolor{blue}{\textbf{Section}}~\ref{sec:4}: under the dense worst-case model, the proposed LP-based formulation and the classical QP-based steepest descent formulation have the same asymptotic complexity order, but the proposed formulation introduces more explicit constraints. Therefore, the classical steepest descent method may have a smaller per-iteration computational cost. However, per-iteration speed alone does not determine overall efficiency. What matters in practice is the computational effort required to reach high-quality objective values under a reasonable computational budget.

\par In the Eastern Massachusetts network, the proposed min--max programming method demonstrates stronger practical convergence behavior than the classical steepest descent algorithm in both Group I and Group II. In Group I, the proposed method reaches \((27807.8324, 29926.3683)\) for the UE and SO objectives within only 1000 iterations and 217.5921 seconds, while the classical steepest descent method remains worse in both objectives even after 5000 iterations and 797.0288 seconds. This indicates that the proposed direction-selection mechanism yields substantially faster and more balanced progress.

\par In Group II, the proposed method reaches a high-quality solution within 1500 iterations and 327.5792 seconds. Under the same 1500-iteration budget, the classical steepest descent method remains substantially worse in both objectives. Even when the classical method is allowed to run for 3000 iterations, which is already twice the iteration count of the proposed method, it still remains worse in both UE and SO objectives. At 5000 iterations, the classical method obtains a slightly lower SO objective, but this improvement is achieved at the expense of a higher UE objective. Therefore, the 5000-iteration steepest descent solution does not dominate the solution obtained by the proposed method. From a practical convergence perspective, the proposed min--max programming method reaches a stable and balanced solution much earlier, whereas the classical steepest descent method requires a substantially larger computational budget to move toward a different trade-off region.

\par In the Anaheim network, the proposed min--max programming method outperforms the classical steepest descent algorithm in both Group I and Group II. As shown in \textcolor{blue}{\textbf{Table}}~\ref{table:3} and \textcolor{blue}{\textbf{Table}}~\ref{table:4}, the proposed method reaches better UE and SO objective values than the classical steepest descent benchmark in both groups. In Group I, the proposed method obtains the best objective values within 2500 iterations and 907.0508 seconds, whereas the classical steepest descent method remains worse even after 20000 iterations and 4669.4956 seconds. In Group II, the proposed method again achieves the best objective values within 2500--3000 iterations, while the classical steepest descent method remains worse after 20000 iterations. In both Anaheim groups, the proposed method uses less than one fifth of the CPU time required by the longest steepest descent run while producing better objective values.

\par Overall, the numerical examples show that the proposed min--max programming method may have a higher per-iteration cost, but it generates more effective common descent directions and reaches high-quality solutions in substantially fewer iterations. This behavior is consistent with the motivation of the proposed direction-finding subproblem: by minimizing the worst directional derivative, the method encourages balanced progress across UE and SO objectives. Consequently, the proposed method provides more reliable practical convergence for large-scale sparse multi-criteria traffic assignment problems, especially when both individual-equilibrium and system-efficiency criteria must be considered.

\begin{table}[ht]
\centering
\caption{Computational efficiency and quality comparison of algorithms across different problem scales: Group I} 
\label{table:3}
\resizebox{\textwidth}{!}{
\begin{tabular}{ccccccc}
\thicktoprule
Network & Index & Algorithm & Iterations & Total CPU Times (s) & UE Obj & SO Obj \\ \hline
EMN & E-1 & Min-Max programming & 1000 & 217.5921 & \textbf{27807.8324} & \textbf{29926.3683} \\
EMN & E-2 & Classical steepest descent  & 1000 & 168.4573 & 29404.6583 & 31450.1727 \\
EMN & E-3 & Classical steepest descent  & 3000 &  497.7244 & 28265.9040 & 30302.8603\\
EMN & E-4 & Classical steepest descent  & 5000 &  797.0288 & 28036.4463 & 30073.4354\\
AN & A-1 & Min-Max programming & 2000 & 689.3975 & \textbf{1369085.7524} & \textbf{1485687.2590} \\
AN & N/A & Min-Max programming & 2500 & 907.0508 & \textbf{1369085.7518} & \textbf{1485687.2586} \\
AN & A-2 & Classical steepest descent  & 2000 & 461.7928 & 1369309.7863 & 1486824.7749 \\
AN & A-3 & Classical steepest descent  & 5000 & 1181.1608 & 1369130.9599 & 1485866.0613 \\
AN & A-4 & Classical steepest descent  & 10000 & 2324.9328 & 1369097.2165 & 1485717.3764 \\
AN & A-5 & Classical steepest descent  & 20000 & 4669.4956 & 1369091.2271 & 1485711.4231 \\
\thickbottomrule
\end{tabular}
}
\end{table}

\begin{table}[htbp]
\centering
\caption{Computational efficiency and quality comparison of algorithms across different problem scales: Group II} 
\label{table:4}
\resizebox{\textwidth}{!}{
\begin{tabular}{ccccccc}
\thicktoprule
Network & Index & Algorithm & Iterations & Total CPU Times (s) & UE Obj & SO Obj \\ \hline
EMN & E-5 & Min-Max programming & 1500 & 327.5792 & \textbf{27783.4166} & 29883.6188 \\
EMN & E-6 & Classical steepest descent  & 1500 & 247.0591 & 30210.3390 & 31911.2084 \\
EMN & E-7 & Classical steepest descent  & 3000 &  497.0678 & 28636.1906 & 30314.5530\\
EMN & E-8 & Classical steepest descent  & 5000 &  784.1923 & 28145.9007 & \textbf{29824.2736}\\
% EMN & new & Classical steepest descent  & 8000 &  1350.9642 & 27995.3012 & 29673.6763\\
% EMN & E-8 & Classical steepest descent  & 10000 &  1676.7893 & 27945.6844 & \textbf{29624.0602}\\
AN & N/A & Min-Max programming & 2000 & 621.9195 & 1369647.1273 & 1486750.4309 \\
AN & A-6 & Min-Max programming & 2500 & 810.7932 & \textbf{1369085.7386} & \textbf{1485695.5007} \\
AN & N/A & Min-Max programming & 3000 & 1007.9466 & \textbf{1369085.7380} & \textbf{1485695.5002} \\
% AN & A-5 & Classical steepest descent  & 1500 & 344.9716 & 1372109.8239 & 1489894.6532 \\
% AN & A-6 & Classical steepest descent  & 2000 & 460.6970 & 1370564.0043 & 1487848.8282 \\
AN & A-7 & Classical steepest descent  & 2500 & 580.9195 & 1369669.8388 & 1486800.2573 \\
AN & A-8 & Classical steepest descent  & 5000 & 1167.7620 & 1369129.9851 & 1485865.2745 \\
AN & A-9 & Classical steepest descent  & 10000 & 2319.9209 & 1369110.9667 & 1485777.6782 \\
AN & A-10 & Classical steepest descent  & 20000 & 4623.1413 & 1369102.9262 & 1485752.7798 \\
\thickbottomrule
\end{tabular}
}
\end{table}

\section{Conclusions and Discussions}

\par This paper presents a min--max multi-gradient descent framework for multi-objective transportation problems. The proposed method identifies a common descent direction by minimizing the maximum directional derivative among all objectives, thereby explicitly emphasizing the least-improved objective at each iteration. Under convex problem-specific constraints, the resulting direction is feasible and provides a balanced first-order improvement whenever a common descent direction exists. Theoretical results establish the existence of the proposed direction and characterize its connection to Pareto criticality.

\par A distinguishing feature of the proposed framework is that, when the problem-specific constraints are linear, the direction-finding subproblem reduces to a linear program. This contrasts with classical multi-objective gradient descent algorithms, such as the steepest descent formulation in \eqref{eq:1} and the minimum-norm formulation in \eqref{eq:2}, which rely on quadratic programming subproblems. The linear-programming formulation does not provide an inherent dense worst-case complexity advantage over the classical QP-based formulations; under a standard dense interior-point complexity model, both formulations have the same asymptotic order. Nevertheless, the proposed formulation directly optimizes the worst directional derivative and therefore provides a different and more balanced direction-selection mechanism. In sparse transportation networks, its practical performance depends on the sparsity pattern, preprocessing, and linear-algebra structure of the resulting Newton systems.

\par To illustrate the practical value of the proposed framework, we apply it to two transportation-related problems: training a physics-informed machine-learning-based car-following model and solving a multi-criteria traffic assignment problem. In the PIML car-following application, the proposed method treats the data-fitting loss and the physics-based loss as two separate objectives and avoids manually fixing scalarization weights. The numerical results show that it achieves the best final predictive performance among all compared training strategies, reducing the average RMSE by approximately \(5.81\%\) relative to the baseline while also reducing variability across random seeds. In the multi-criteria traffic assignment application, which balances user equilibrium and system optimum objectives, the proposed method exhibits better convergence behavior than the classical steepest descent formulation. On the tested transportation networks, it reaches high-quality solutions with substantially fewer iterations and lower total CPU time, demonstrating the practical advantage of the min--max direction-selection mechanism in large-scale sparse settings.

\par Beyond the current convex and linearly constrained framework, future research could explore extensions of the proposed min-max MGDA idea to nonconvex constrained multi-objective optimization problems. In the present study, the main constrained transportation application involves linear constraints, so feasibility can be checked directly and feasible steps can be selected reliably. For more general differentiable nonlinear constraints, one possible extension is to replace the hard feasibility constraint in the direction-finding subproblem with a local linearized approximation. For example, at a current feasible point \(\boldsymbol{\theta}\), one may consider the following linearized direction-finding subproblem:
\begin{subequations}\label{eq:109}
\begin{align}
  &w
    =
    \min_{\mathbf{d}}\;\max_{t=1,2,\dots,T}
      \nabla\mathcal{F}_{t}(\boldsymbol{\theta})^{\top}\mathbf{d}
    \label{eq:109a}\\
  &\text{s.t.}\quad
    f_{i}(\boldsymbol{\theta})
    +
    \nabla f_{i}(\boldsymbol{\theta})^{\top}\mathbf{d}
    \le 0,
    \quad i=1,2,\dots,m
    \label{eq:109b}\\
  &\quad
    \nabla \mathcal{F}_{t}(\boldsymbol{\theta})^{\top}\mathbf{d} \le 0,
    \quad t=1,2,\dots,T
    \label{eq:109c}\\
  &\quad
    -\mathbf{1} \le \mathbf{d} \le \mathbf{1}
    \label{eq:109d}
\end{align}
\end{subequations}
This formulation should be interpreted only as a local approximation of the nonlinear feasible region. Since the original nonlinear constraints may still be violated after taking a finite step, the linearized subproblem must be combined with an acceptance or rejection mechanism. Specifically, after obtaining a candidate direction \(\mathbf d\), one can choose a trial step size \(h>0\) and check the original constraints
\begin{equation}\label{eq:110}
f_i(\boldsymbol{\theta}+h\mathbf d)\le 0,
\qquad i=1,\ldots,m
\end{equation}
If all hard constraints are satisfied, the update is accepted; otherwise, the step size is reduced, for example by setting \(h\leftarrow \rho h\) with \(\rho\in(0,1)\), until feasibility is restored or a minimum step-size tolerance is reached. Developing rigorous convergence guarantees for such linearization-and-acceptance schemes under nonconvex constraints is an important direction for future work.

\section{CRediT} 
\textbf{Yuan-Zheng Lei}: Conceptualization, Methodology, Writing - original draft. \textbf{Yaobang Gong}: Conceptualization, Writing, and Experiment design.  \textbf{Xianfeng Terry Yang}: Conceptualization, Methodology, Supervision, and Writing.

\section{Acknowledgement} 
This research is supported by the award \textit{CAREER: Physics Regularized Machine Learning Theory: Modeling Stochastic Traffic Flow Patterns for Smart Mobility Systems (\# 2234289)}, which is funded by the National Science Foundation.

\section{Appendix} \label{8}
\subsection{The convergence curves of the experiments reported in \textcolor{blue}{\textbf{Tables}}~\ref{table:2} and \ref{table:3}}
\begin{figure}[htbp]
    \centering
    \includegraphics[width=0.85\textwidth]{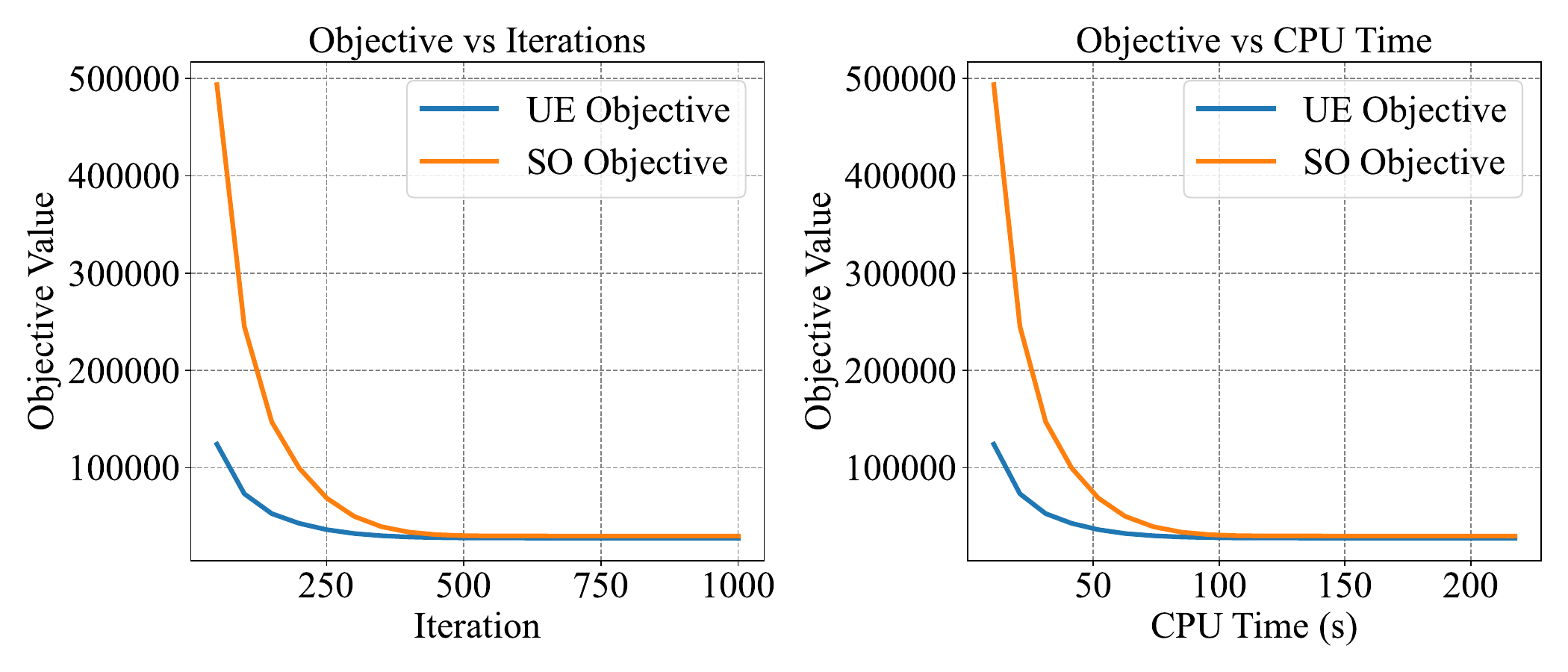}
    \caption{Comparison of UE and SO objective values over iterations using the min-max programming on the Eastern Massachusetts network (E-1).}
    \label{fig:9}
\end{figure}

\begin{figure}[htbp]
    \centering
    \includegraphics[width=0.85\textwidth]{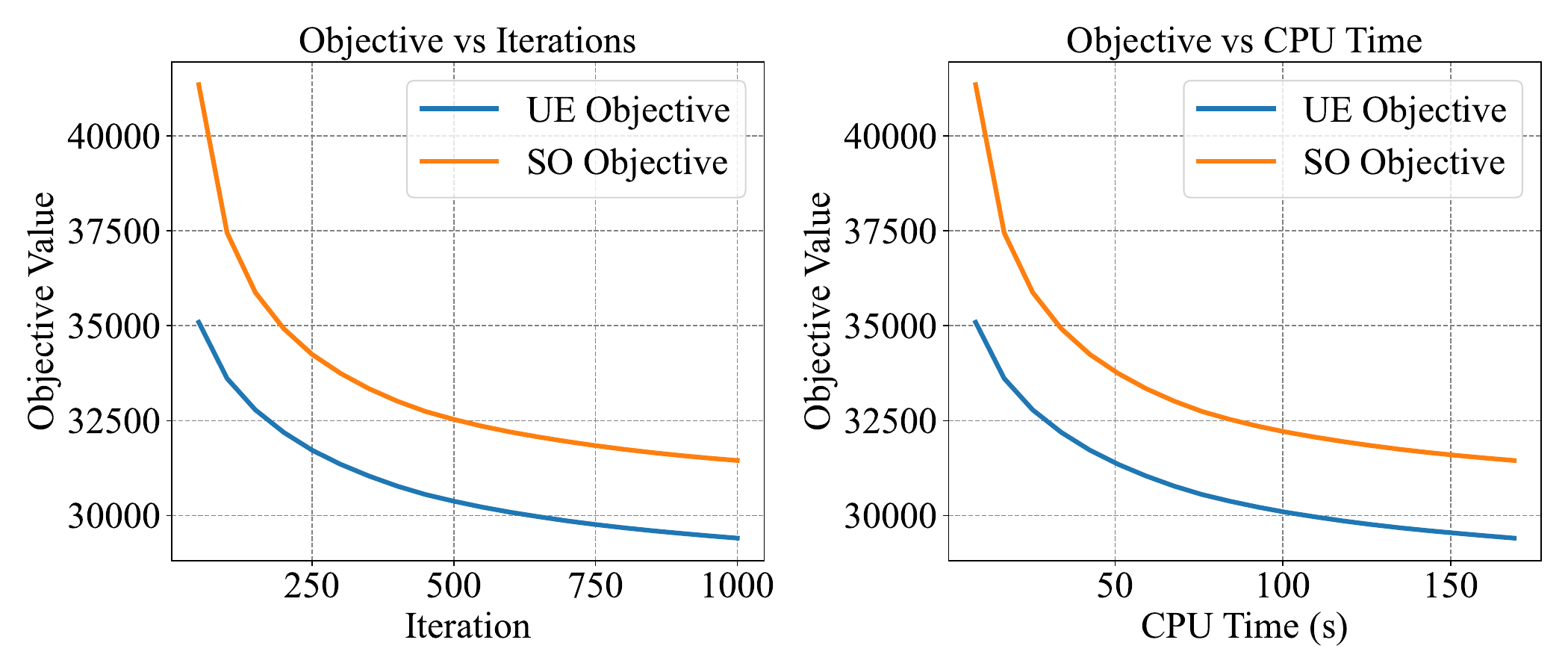}
    \caption{Comparison of UE and SO objective values over iterations using the classical steepest descent method \cite{fliege2000steepest} on the Eastern Massachusetts network (E-2).}
    \label{fig:10}
\end{figure}

\begin{figure}[htbp]
    \centering
    \includegraphics[width=0.85\textwidth]{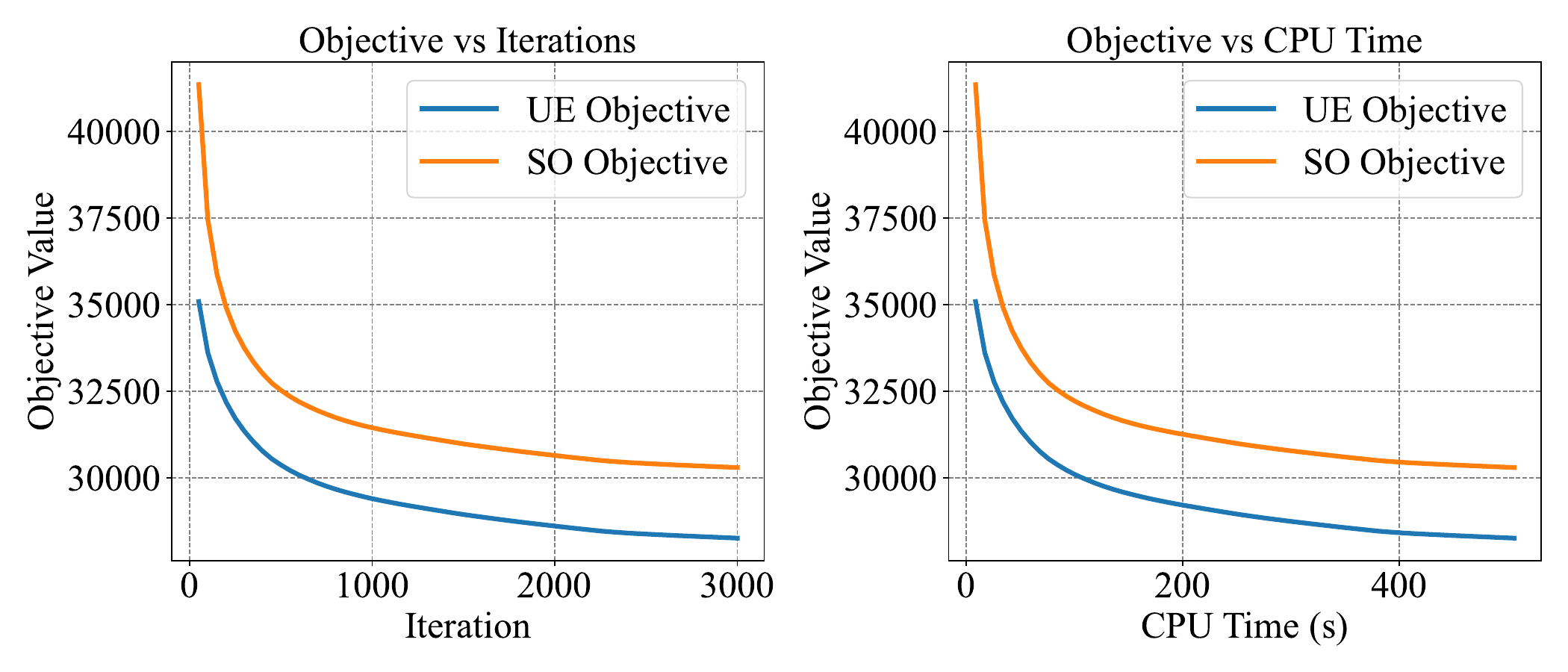}
    \caption{Comparison of UE and SO objective values over iterations using the classical steepest descent method \cite{fliege2000steepest} on the Eastern Massachusetts network (E-3).}
    \label{fig:11}
\end{figure}

\begin{figure}[htbp]
    \centering
    \includegraphics[width=0.85\textwidth]{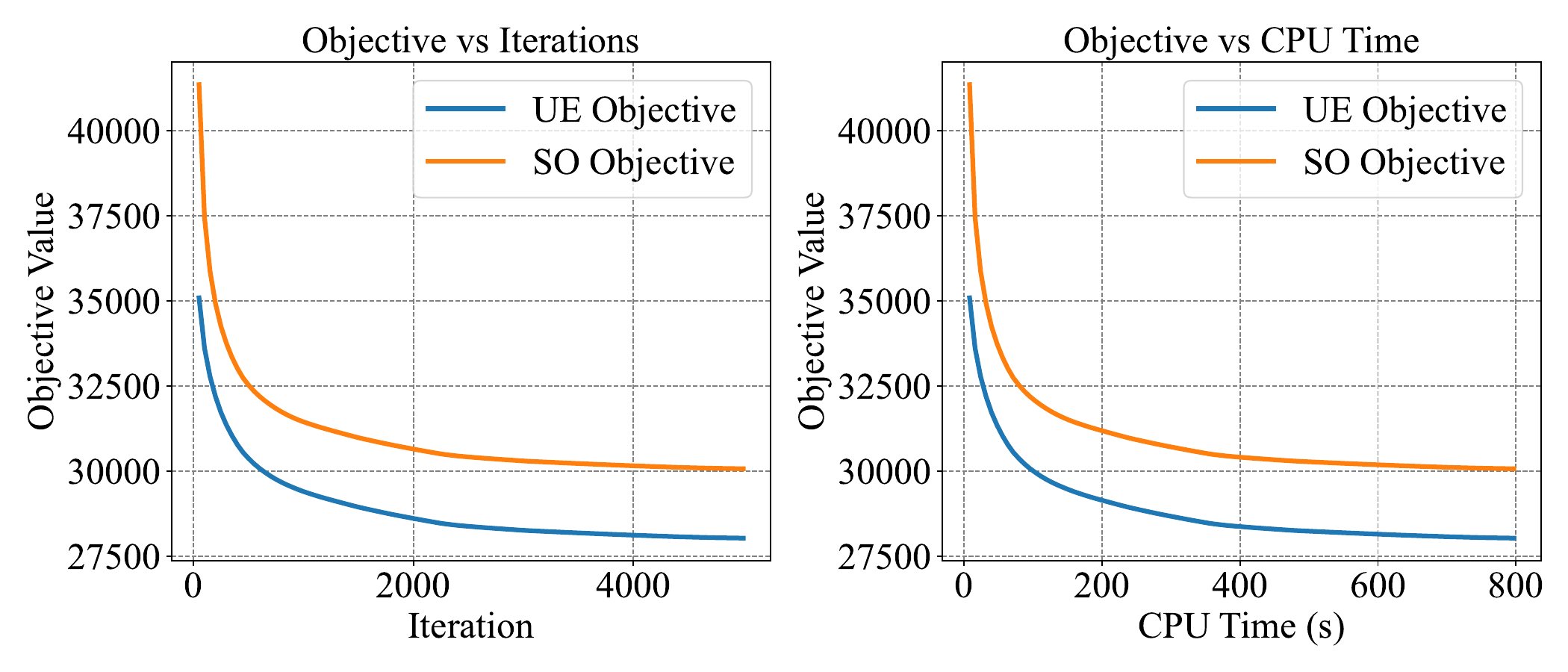}
    \caption{Comparison of UE and SO objective values over iterations using the classical steepest descent method \cite{fliege2000steepest} on the Eastern Massachusetts network (E-4).}
    \label{fig:12}
\end{figure}

\begin{figure}[htbp]
    \centering
    \includegraphics[width=0.85\textwidth]{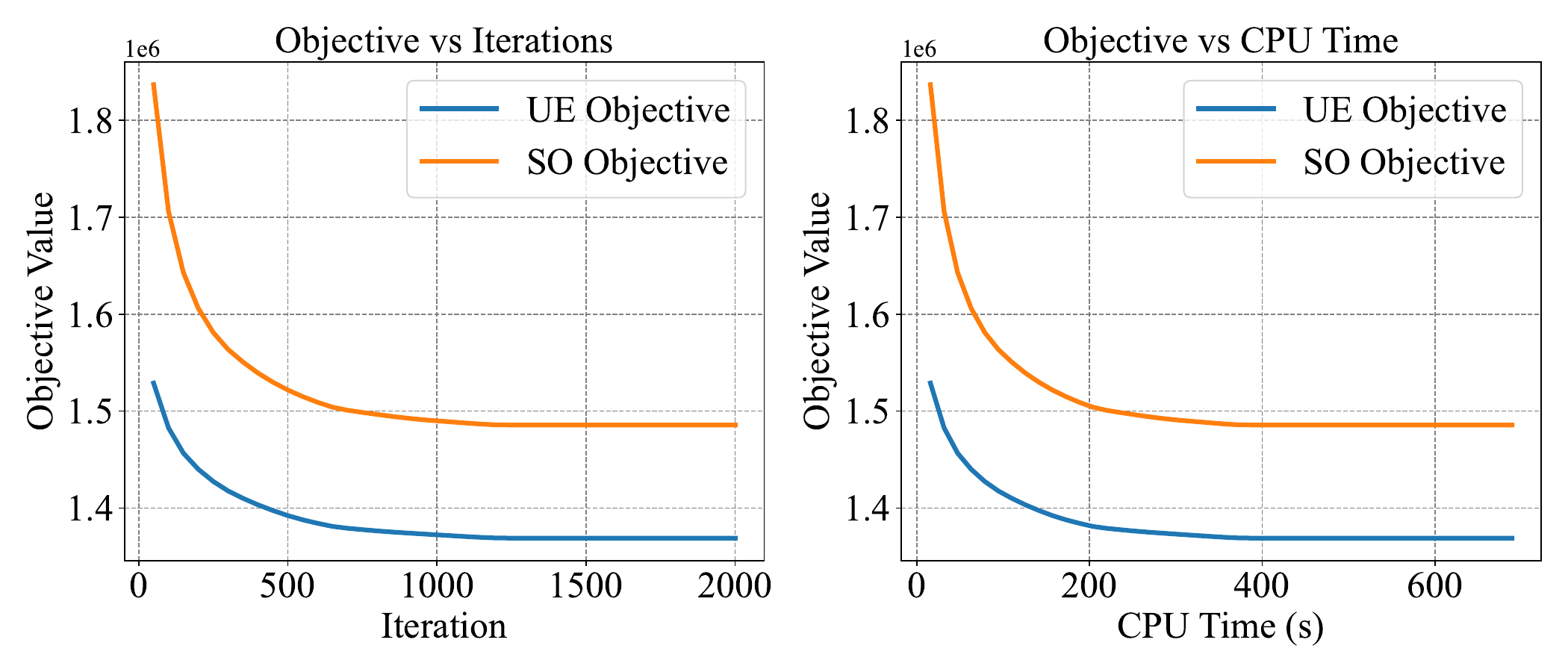}
    \caption{Comparison of UE and SO objective values over iterations using the min-max programming on the Anaheim network (A-1).}
    \label{fig:13}
\end{figure}

\begin{figure}[htbp]
    \centering
    \includegraphics[width=0.85\textwidth]{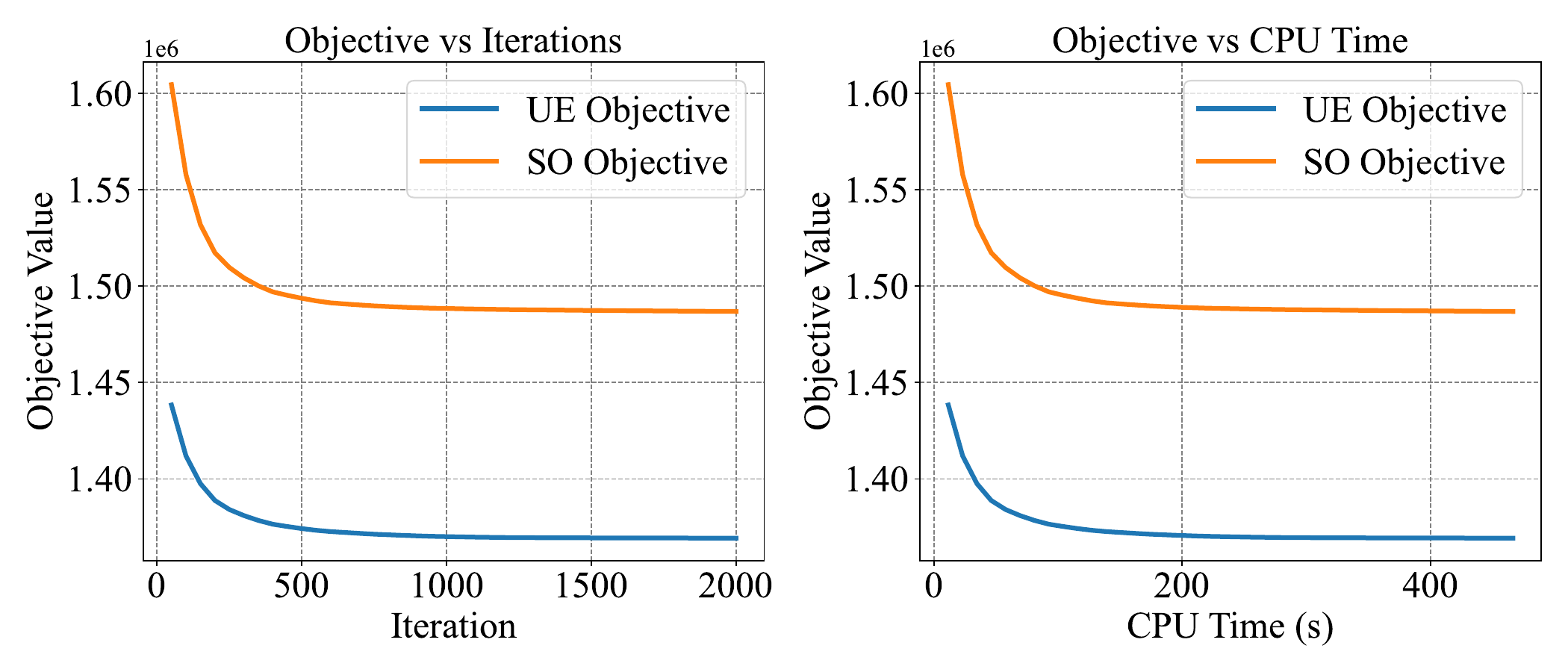}
    \caption{Comparison of UE and SO objective values over iterations using the classical steepest descent method \cite{fliege2000steepest} on the Anaheim network (A-2).}
    \label{fig:14}
\end{figure}

\begin{figure}[htbp]
    \centering
    \includegraphics[width=0.85\textwidth]{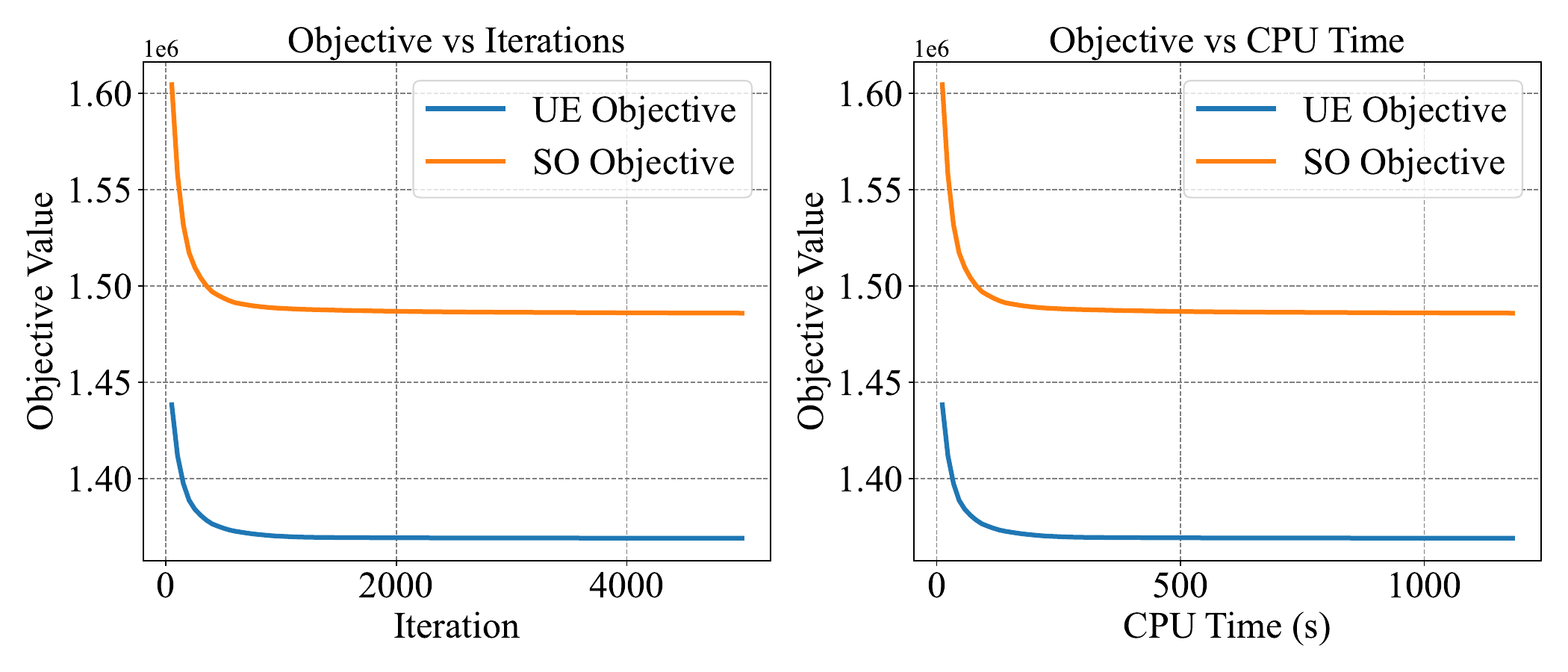}
    \caption{Comparison of UE and SO objective values over iterations using the classical steepest descent method \cite{fliege2000steepest} on the Anaheim network (A-3).}
    \label{fig:15}
\end{figure}

\begin{figure}[htbp]
    \centering
    \includegraphics[width=0.85\textwidth]{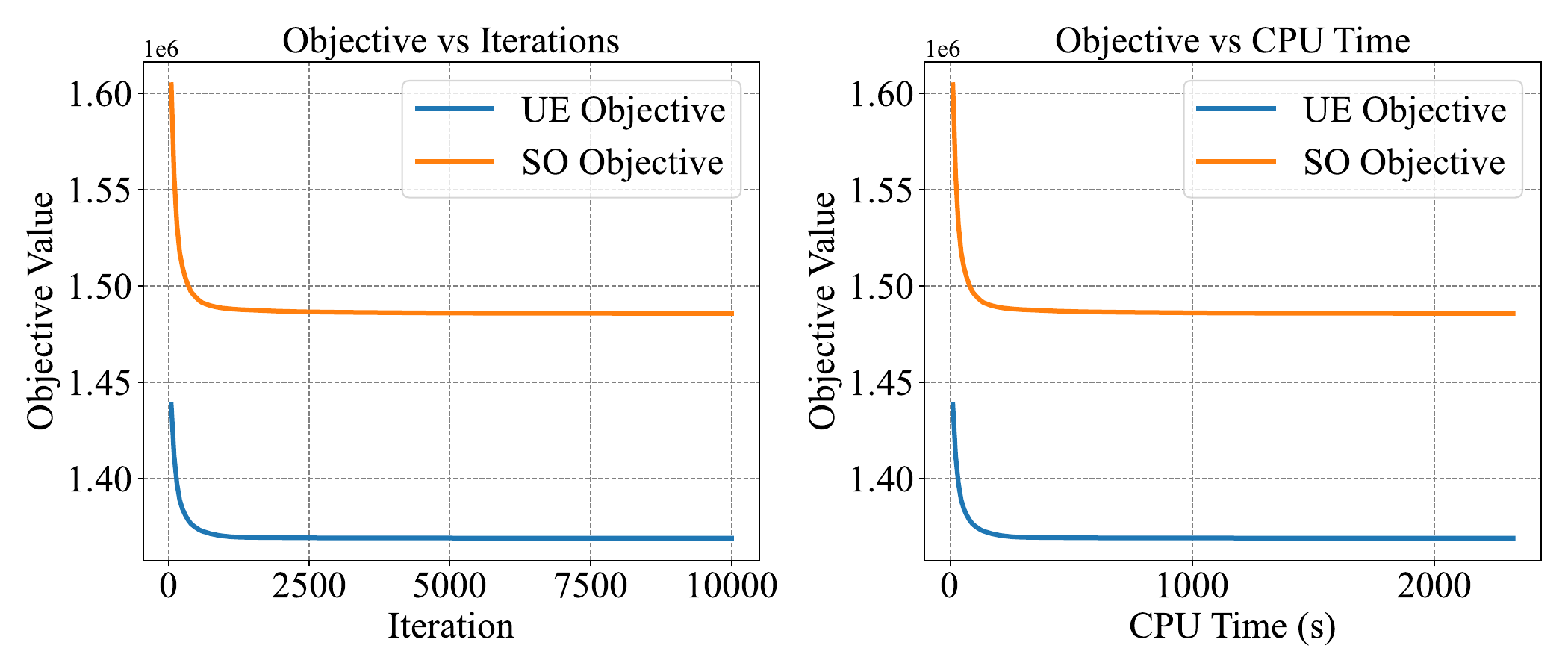}
    \caption{Comparison of UE and SO objective values over iterations using the classical steepest descent method \cite{fliege2000steepest} on the Anaheim network (A-4).}
    \label{fig:16}
\end{figure}

\begin{figure}[htbp]
    \centering
    \includegraphics[width=0.85\textwidth]{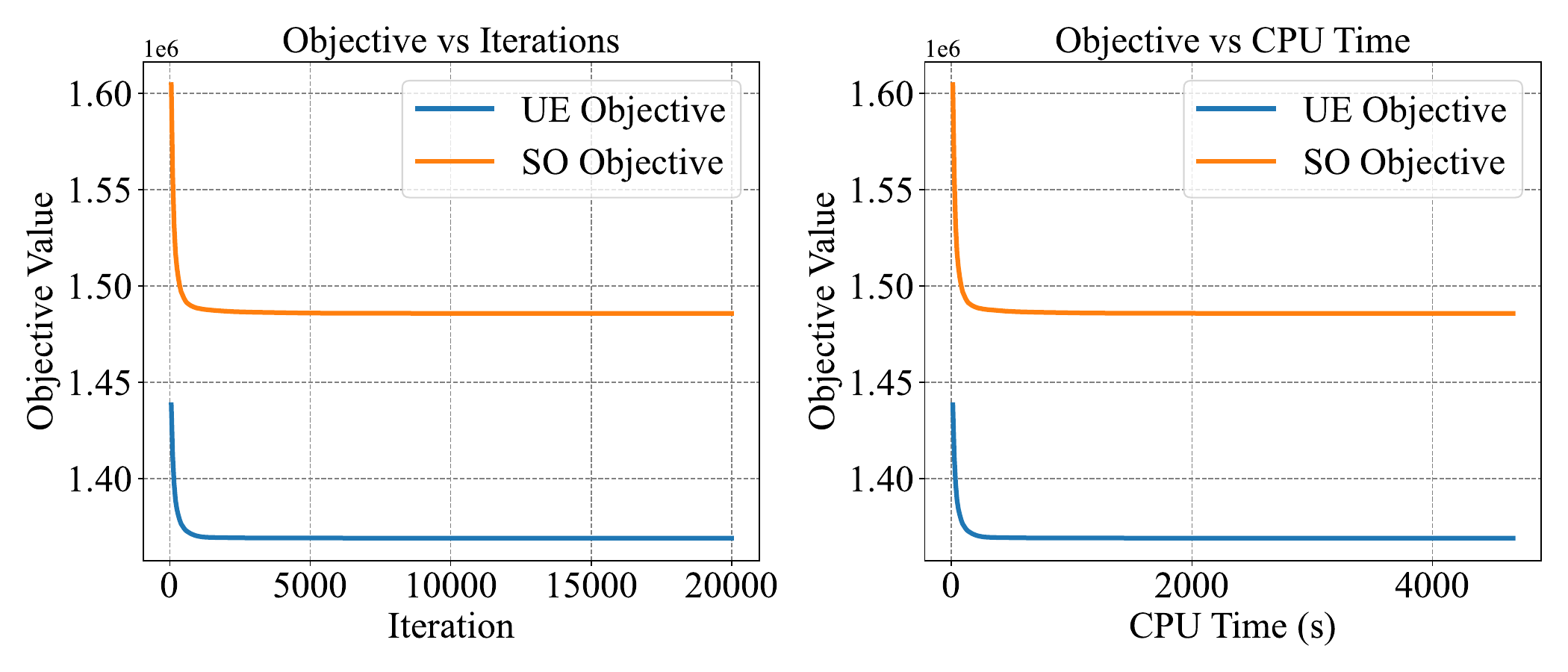}
    \caption{Comparison of UE and SO objective values over iterations using the classical steepest descent method \cite{fliege2000steepest} on the Anaheim network (A-5).}
    \label{fig:17}
\end{figure}

\begin{figure}[htbp]
    \centering
    \includegraphics[width=0.85\textwidth]{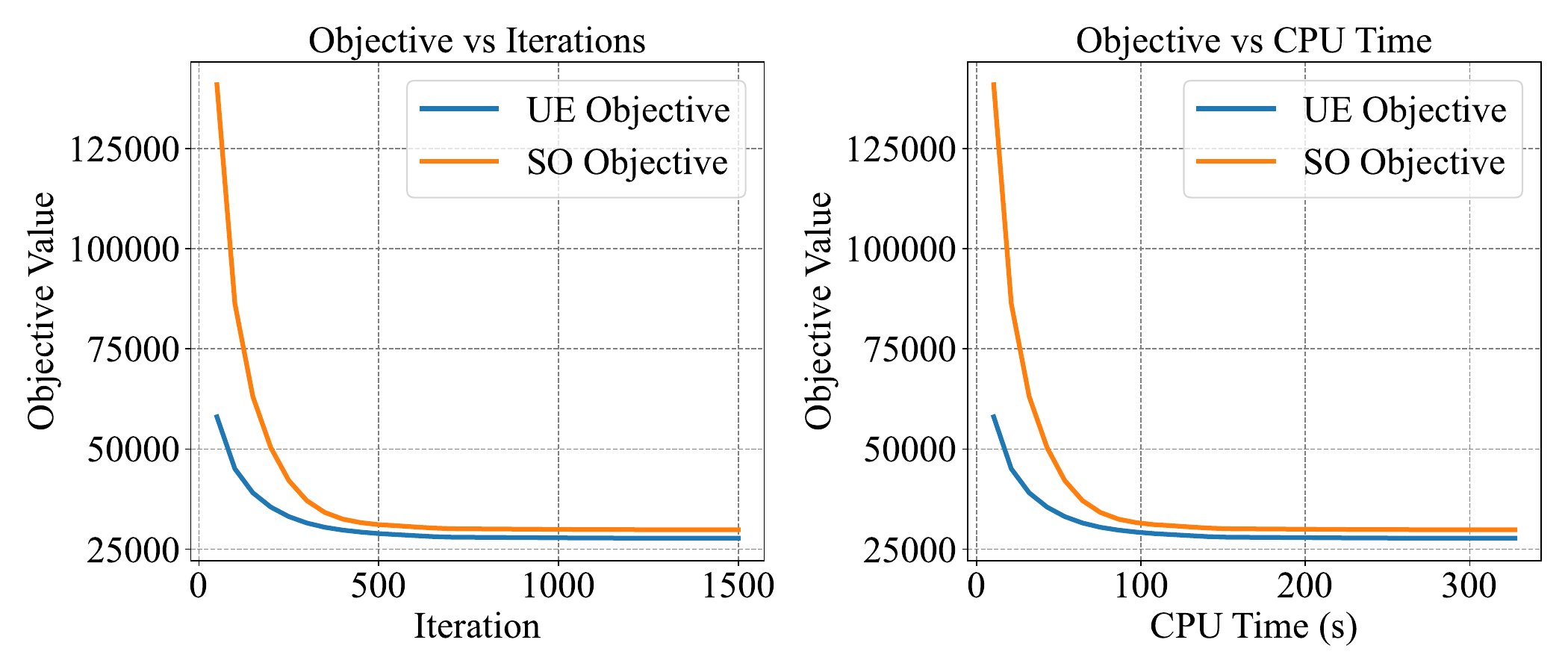}
    \caption{Comparison of UE and SO objective values over iterations using the min-max programming on the Eastern Massachusetts network (E-5).}
    \label{fig:18}
\end{figure}

\begin{figure}[htbp]
    \centering
    \includegraphics[width=0.85\textwidth]{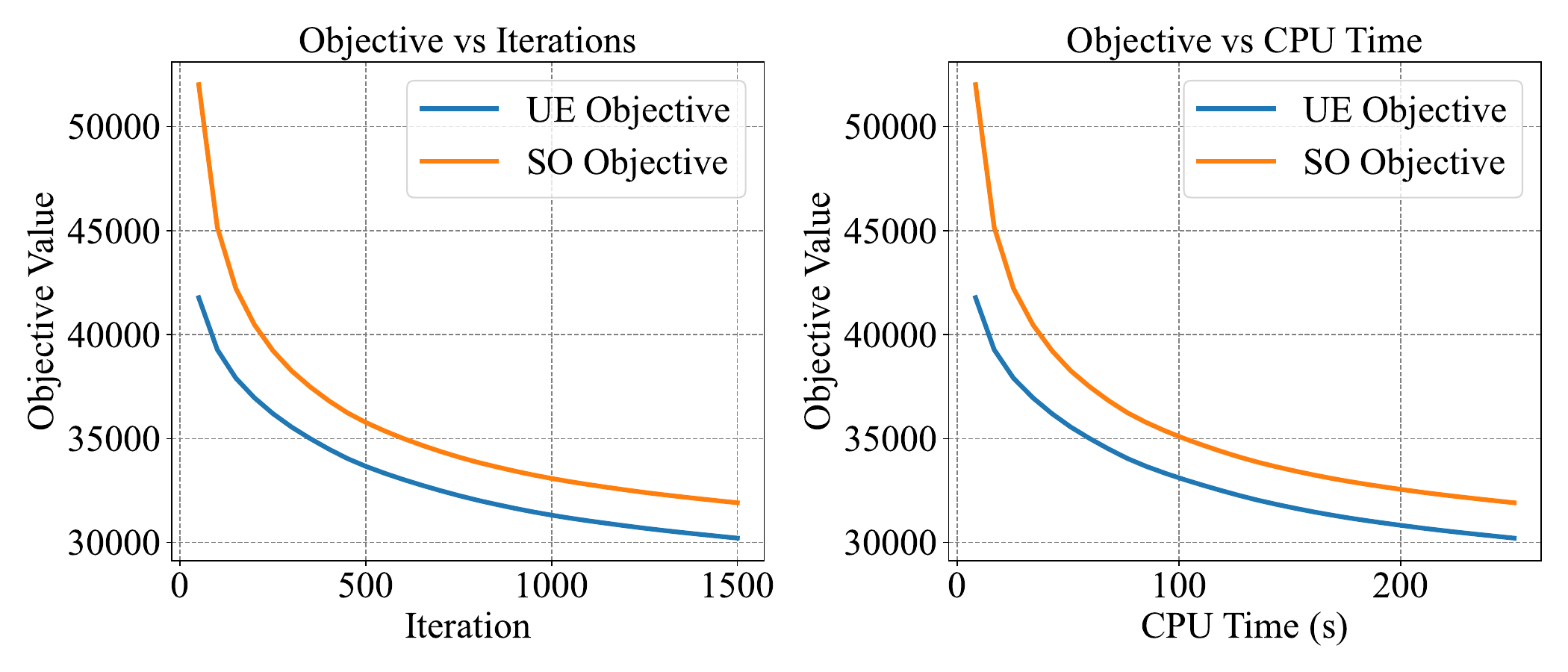}
    \caption{Comparison of UE and SO objective values over iterations using the classical steepest descent method \cite{fliege2000steepest} on the Eastern Massachusetts network (E-6).}
    \label{fig:19}
\end{figure}

\begin{figure}[htbp]
    \centering
    \includegraphics[width=0.85\textwidth]{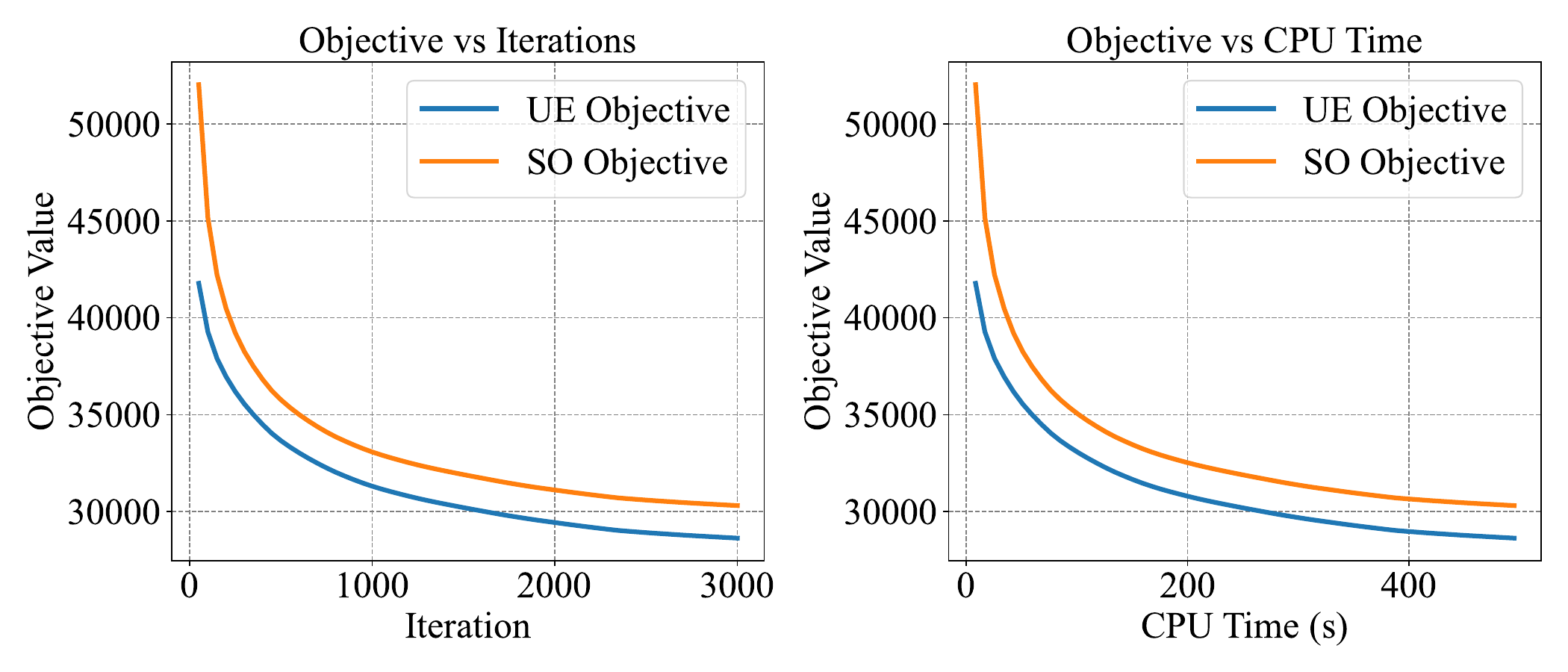}
    \caption{Comparison of UE and SO objective values over iterations using the classical steepest descent method \cite{fliege2000steepest} on the Eastern Massachusetts network (E-7).}
    \label{fig:20}
\end{figure}

\begin{figure}[htbp]
    \centering
    \includegraphics[width=0.85\textwidth]{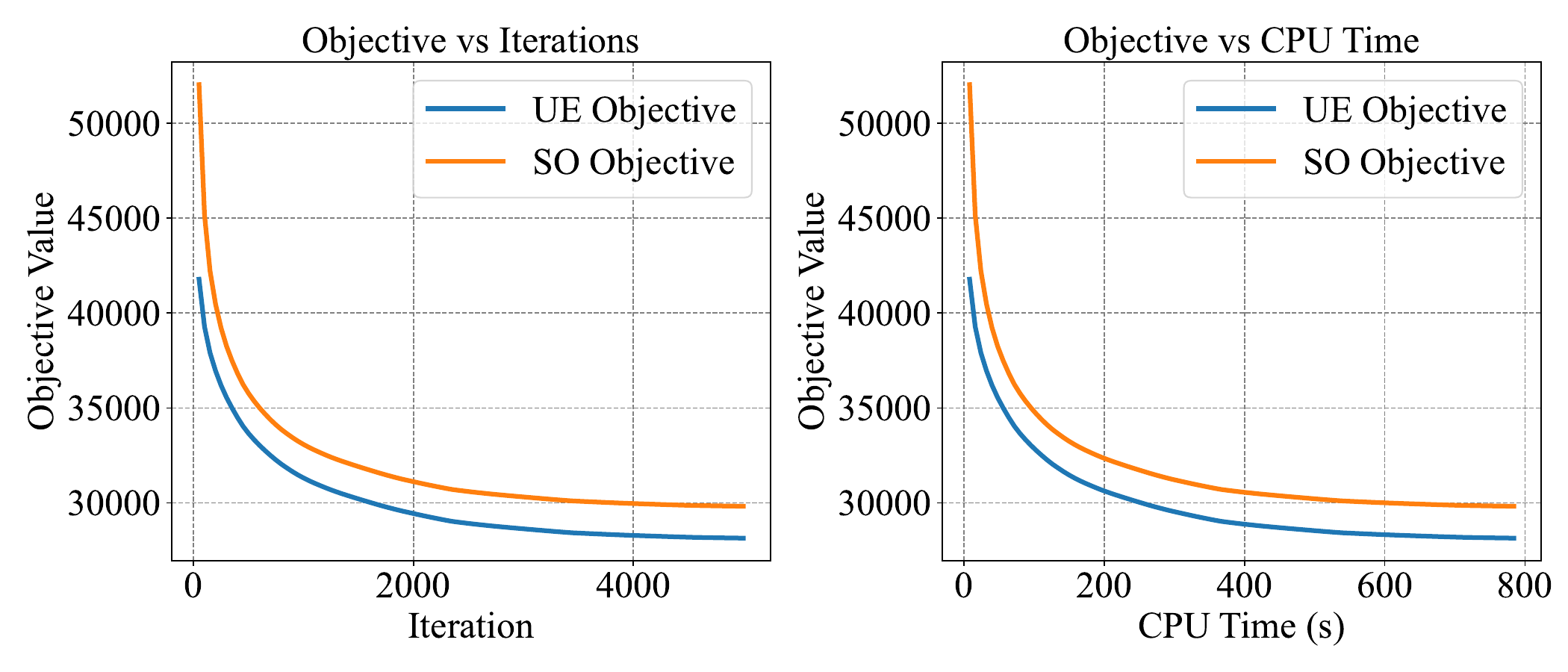}
    \caption{Comparison of UE and SO objective values over iterations using the classical steepest descent method \cite{fliege2000steepest} on the Eastern Massachusetts network (E-8).}
    \label{fig:21}
\end{figure}

\begin{figure}[htbp]
    \centering
    \includegraphics[width=0.85\textwidth]{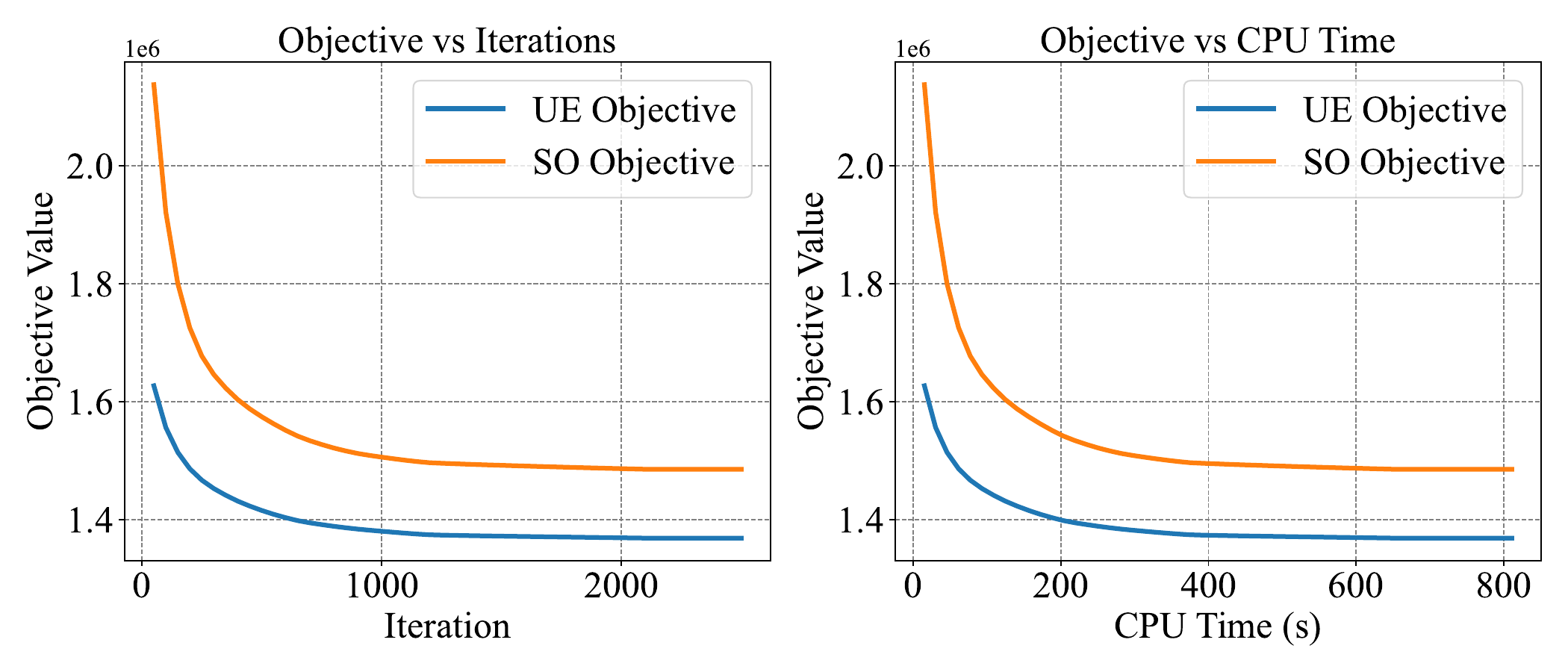}
    \caption{Comparison of UE and SO objective values over iterations using the min-max programming on the Anaheim network (A-6).}
    \label{fig:22}
\end{figure}

\begin{figure}[htbp]
    \centering
    \includegraphics[width=0.85\textwidth]{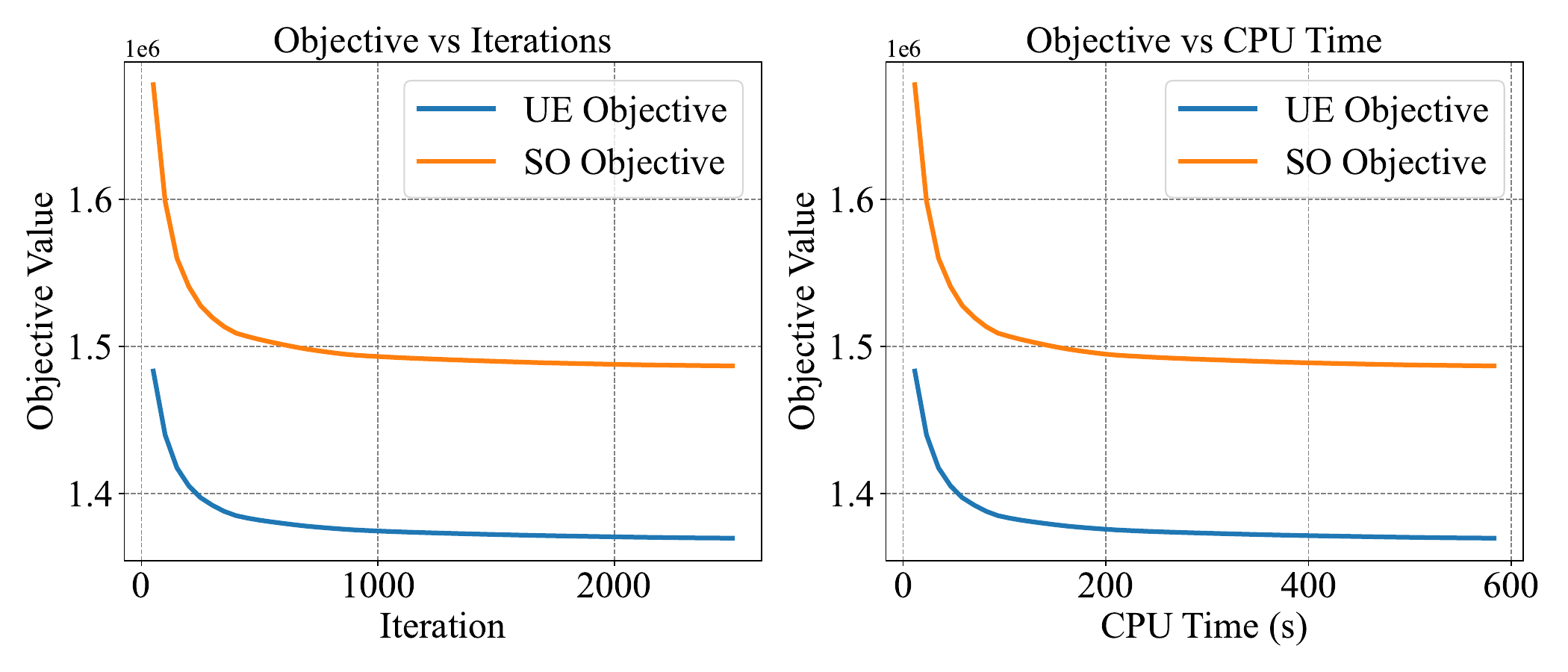}
    \caption{Comparison of UE and SO objective values over iterations using the classical steepest descent method \cite{fliege2000steepest} on the Anaheim network (A-7).}
    \label{fig:23}
\end{figure}

\begin{figure}[htbp]
    \centering
    \includegraphics[width=0.85\textwidth]{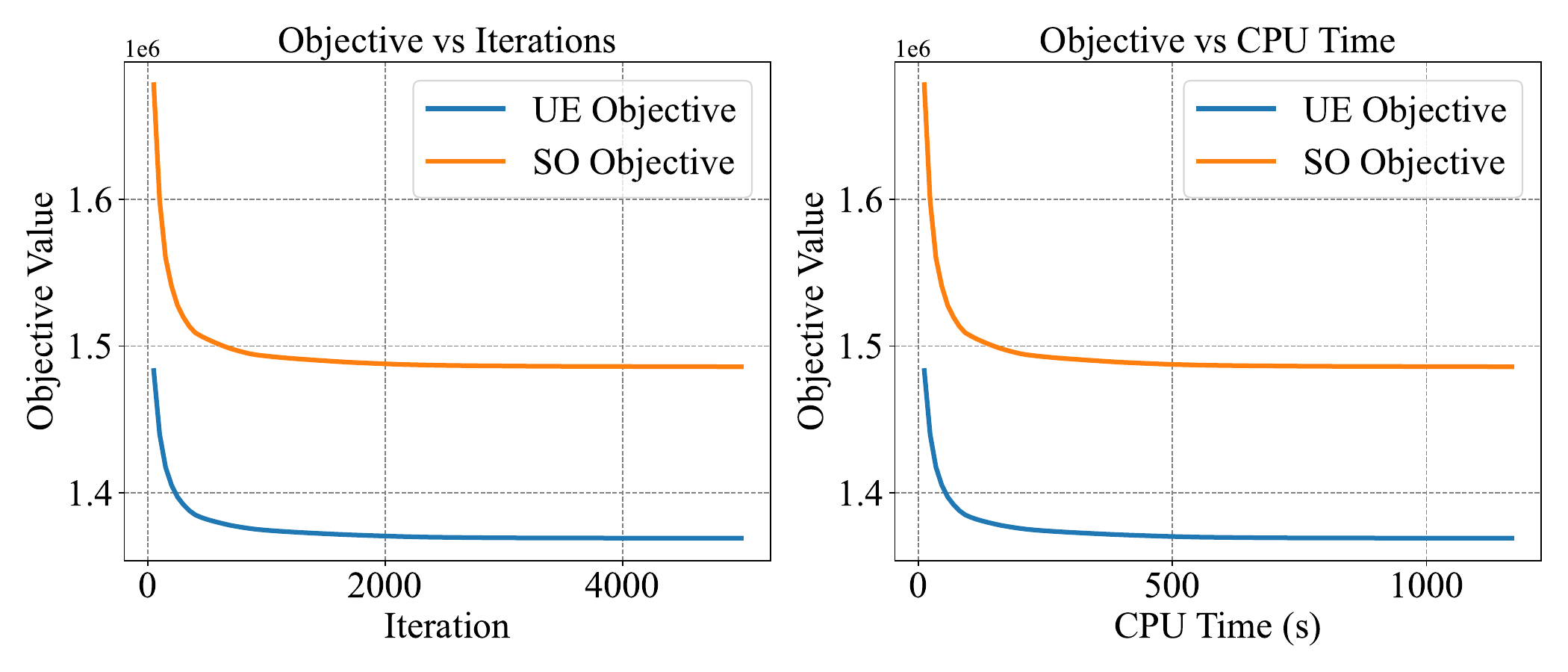}
    \caption{Comparison of UE and SO objective values over iterations using the classical steepest descent method \cite{fliege2000steepest} on the Anaheim network (A-8).}
    \label{fig:24}
\end{figure}

\begin{figure}[htbp]
    \centering
    \includegraphics[width=0.85\textwidth]{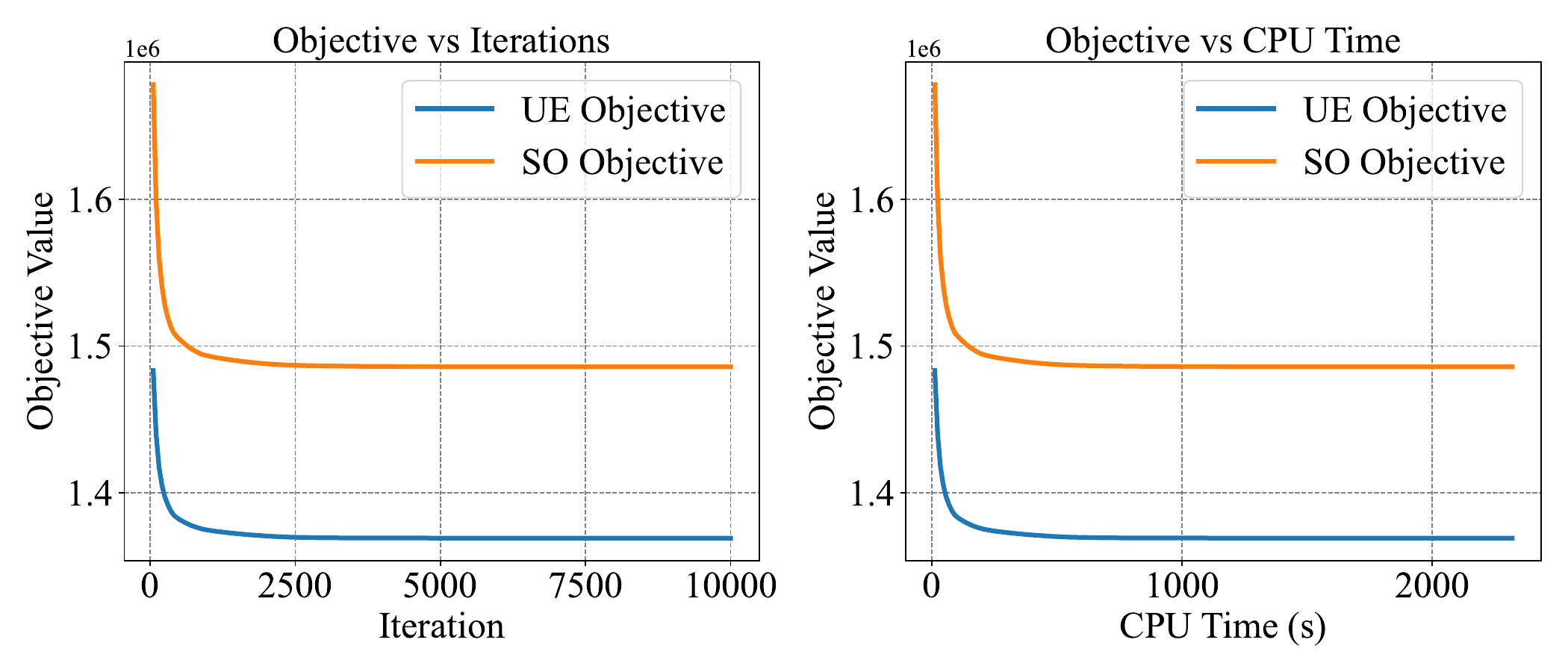}
    \caption{Comparison of UE and SO objective values over iterations using the classical steepest descent method \cite{fliege2000steepest} on the Anaheim network (A-9).}
    \label{fig:25}
\end{figure}

\begin{figure}[htbp]
    \centering
    \includegraphics[width=0.85\textwidth]{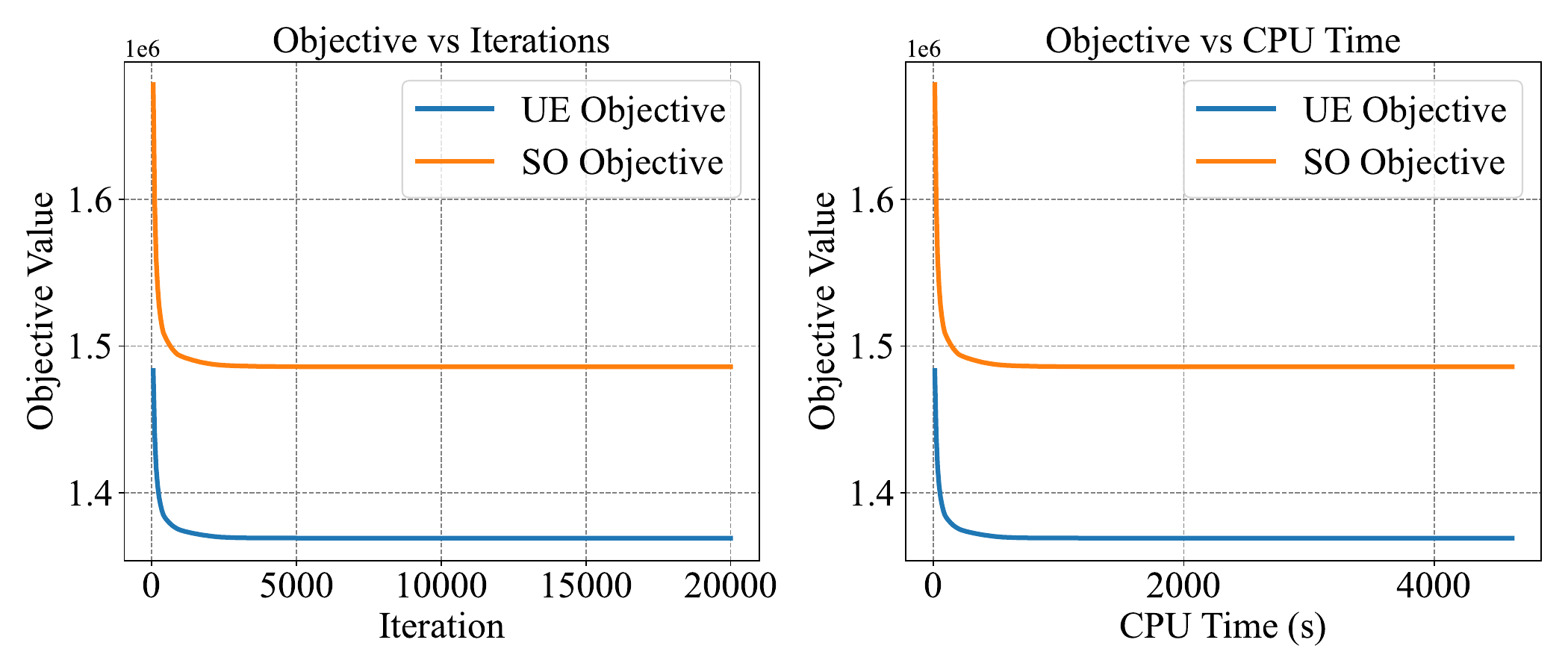}
    \caption{Comparison of UE and SO objective values over iterations using the classical steepest descent method \cite{fliege2000steepest} on the Anaheim network (A-10).}
    \label{fig:26}
\end{figure}

\newpage
\subsection{\textcolor{blue}{\textbf{Theorem}} \ref{theorem:3} and its proof}
\begin{theorem} \label{theorem:3}
Let
\begin{equation} \label{eq:111}
g_t \coloneqq \nabla F_t(\theta)\in\mathbb{R}^n\qquad t=1,\dots,T
\end{equation}
Consider the quadratic program
\begin{equation}
\begin{aligned}
\qquad
\min_{d\in\mathbb{R}^n,\ \gamma\in\mathbb{R}}
\quad & \gamma + \frac12\|d\|^2 \\
\textnormal{s.t.}\quad
& g_t^\top d \le \gamma\qquad t=1,\dots,T
\end{aligned}
\label{eq:112}
\end{equation}
Then its Lagrange dual is
\begin{equation}
\begin{aligned}
\qquad
\max_{\mu\in\mathbb{R}^T}
\quad & -\frac12\left\|\sum_{t=1}^T \mu_t g_t\right\|^2 \\
\textnormal{s.t.}\quad
& \sum_{t=1}^T \mu_t = 1, \\
& \mu_t \ge 0,\qquad t=1,\dots,T.
\end{aligned}
\label{eq:113}
\end{equation}
Equivalently, \eqref{eq:113} can be written as
\begin{equation}
\begin{aligned}
\min_{\mu\in\mathbb{R}^T}
\quad & \left\|\sum_{t=1}^T \mu_t g_t\right\|^2 \\
\textnormal{s.t.}\quad
& \sum_{t=1}^T \mu_t = 1 \\
& \mu_t \ge 0,\qquad t=1,\dots,T
\end{aligned}
\label{eq:114}
\end{equation}
which is exactly Problem \eqref{eq:2}. Hence, Problem \eqref{eq:2} is an equivalent minimization form of the Lagrange dual of Problem \eqref{eq:1}.

Moreover, if $\mu^*$ solves \eqref{eq:114}, then the primal optimal solution satisfies
\begin{equation} \label{eq;115}
d^* = -\sum_{t=1}^T \mu_t^* g_t
\qquad
\gamma^* = -\left\|\sum_{t=1}^T \mu_t^* g_t\right\|^2
\end{equation}
\end{theorem}

\begin{proof}
The Lagrangian of \eqref{eq:112} is
\begin{equation} \label{eq:116}
\mathcal{L}(d,\gamma,\mu)
=
\gamma+\frac12\|d\|^2
+\sum_{t=1}^T \mu_t\bigl(g_t^\top d-\gamma\bigr)
\end{equation}
where $\mu_t\ge 0$ for all $t=1,\dots,T$. Rearranging terms gives
\begin{equation} \label{eq;117}
\mathcal{L}(d,\gamma,\mu)
=
\frac12\|d\|^2
+\left(\sum_{t=1}^T \mu_t g_t\right)^\top d
+\Bigl(1-\sum_{t=1}^T\mu_t\Bigr)\gamma
\end{equation}

We now compute the dual function
\begin{equation} \label{eq:118}
q(\mu)=\inf_{d,\gamma}\mathcal{L}(d,\gamma,\mu)
\end{equation}

First, if $\sum_{t=1}^T\mu_t\neq 1$, then the coefficient of $\gamma$ is nonzero, and thus
\begin{equation}  \label{eq:119}
\inf_{\gamma\in\mathbb{R}}
\Bigl(1-\sum_{t=1}^T\mu_t\Bigr)\gamma=-\infty
\end{equation}
Hence, for the dual function to be finite, it is necessary that
\begin{equation} \label{eq:120}
\sum_{t=1}^T\mu_t=1
\end{equation}

Under this condition, the Lagrangian reduces to
\begin{equation} \label{eq:121}
\mathcal{L}(d,\gamma,\mu)
=
\frac12\|d\|^2
+\left(\sum_{t=1}^T \mu_t g_t\right)^\top d
\end{equation}
Let
\begin{equation} \label{eq:122}
s(\mu)\coloneqq \sum_{t=1}^T \mu_t g_t
\end{equation}
Then
\begin{equation} \label{eq:123}
\mathcal{L}(d,\gamma,\mu)=\frac12\|d\|^2+s(\mu)^\top d
\end{equation}
Minimizing with respect to $d$ yields
\begin{equation} \label{eq:124}
\nabla_d \mathcal{L}(d,\gamma,\mu)=d+s(\mu)=0
\end{equation}
so the minimizing $d$ is
\begin{equation} \label{eq;125}
d^*(\mu)=-s(\mu)=-\sum_{t=1}^T \mu_t g_t
\end{equation}
Substituting this back gives
\begin{equation} \label{eq:126}
q(\mu)
=
\inf_d\left(\frac12\|d\|^2+s(\mu)^\top d\right)
=
-\frac12\|s(\mu)\|^2
=
-\frac12\left\|\sum_{t=1}^T \mu_t g_t\right\|^2
\end{equation}
Therefore, the dual problem is exactly
\begin{equation} \label{eq:127}
\max_{\mu_t\ge 0,\ \sum_t\mu_t=1}
-\frac12\left\|\sum_{t=1}^T \mu_t g_t\right\|^2
\end{equation}
which is \eqref{eq:113}. Since multiplying the objective by the positive constant $2$ and changing from maximizing a negative quantity to minimizing its opposite do not alter the optimizer, \eqref{eq:114} is equivalent to
\begin{equation}
\min_{\mu_t\ge 0,\ \sum_t\mu_t=1}
\left\|\sum_{t=1}^T \mu_t g_t\right\|^2
\end{equation}
which is \eqref{eq:114}, namely Problem (2).

Finally, since \eqref{eq:113} is a convex program and satisfies Slater's condition
(for example, $d=0$ and any $\gamma>0$ is strictly feasible), strong duality holds. Hence any optimal dual solution $\mu^*$ recovers the primal optimal direction as
\begin{equation} \label{eq:129}
d^*=-\sum_{t=1}^T \mu_t^* g_t
\end{equation}
Moreover, by complementary slackness,
\begin{equation} \label{eq:130}
\sum_{t=1}^T \mu_t^*(g_t^\top d^*-\gamma^*)=0
\end{equation}
Using $\sum_{t=1}^T\mu_t^*=1$, we obtain
\begin{equation} \label{eq:131}
\gamma^*=\sum_{t=1}^T \mu_t^* g_t^\top d^*
=\left(\sum_{t=1}^T \mu_t^* g_t\right)^\top d^*
=-\left\|\sum_{t=1}^T \mu_t^* g_t\right\|^2
\end{equation}
This completes the proof.
\end{proof}

% \section{Model Source Code} \label{10}
% The source code of the proposed algorithm will be made public after the peer-review process.

% Insert cited references

%\newpage
\bibliography{ref}

\end{document}